\documentclass[reqno,11pt]{amsart}

\usepackage{enumerate}
\usepackage{amssymb}
\usepackage{fullpage}

\theoremstyle{plain}
\newtheorem{thm}{Theorem}[section]
\newtheorem{prop}[thm]{Proposition} 
\newtheorem{lem}[thm]{Lemma} 
\newtheorem{cor}[thm]{Corollary} 
\theoremstyle{definition}
\newtheorem{defn}[thm]{Definition} 
\theoremstyle{remark}
\newtheorem{rem}[thm]{\textsl{Remark}} 
\newtheorem{exmp}[thm]{\textsl{Example}} 
\numberwithin{equation}{section}

\DeclareMathAlphabet{\mathsfsl}{OT1}{cmss}{m}{sl}

\newcommand\ad{\operatorname{ad}}
\newcommand\appli[5]{
  \begin{matrix}#1\colon\vphantom{#2#3}\\*\vphantom{#4#5}\end{matrix}
  \begin{matrix}#2\vphantom{#1#3}\\*#4\vphantom{#5}\end{matrix}\;
  \begin{matrix}\longrightarrow\vphantom{#1#2#3}\\*
    \longmapsto\vphantom{#4#5}\end{matrix}\;
  \begin{matrix}#3\vphantom{#1#2}\\*#5\vphantom{#4}\end{matrix}}
\newcommand\barGamma{{\overline{\Gamma}}}
\newcommand\can{\mathit{can}}
\renewcommand\cosh{\operatorname{ch}}
\renewcommand\coth{\operatorname{coth}}
\newcommand\cprime{$'$}
\renewcommand\d{\operatorname{d}}
\newcommand\dlie{\lie{d}}
\newcommand\dell{\partial}
\newcommand\e{\operatorname{e}}
\newcommand\eg{e.g.,\ }
\newcommand\ensemble[2]{\left.\left\{#1\vphantom{#2}\,
    \right\vert\,#2\right\}}
\newcommand\g{\lie{g}}
\newcommand\h{\lie{h}}
\newcommand\hfills{\hspace*{\fill}}
\newcommand\ie{i.e.,\ }
\newcommand\im{\mathrm{i}}
\renewcommand\l{\lie{l}}
\newcommand\lcan{l^\can}
\newcommand\lcanG[1]{l^{\can(#1)}}
\newcommand\lie[1]{\mathfrak{#1}}
\newcommand\m{\lie{m}}
\newcommand\n{\lie{n}}
\newcommand\op{\mathsf{op}}
\newcommand\p{\mathrm{p}}

\newcommand\polhk[1]{\setbox0=\hbox{#1}{\ooalign{\hidewidth
      \lower1.5ex\hbox{`}\hidewidth\crcr\unhbox0}}}
\newcommand\rank{\mathsfsl{rank}}

\newcommand\rAM{r^{\mathit{AM}}}
\renewcommand\sinh{\operatorname{sh}}
\renewcommand\t{\mathsfsl{t\mskip.2\thinmuskip}}
\renewcommand\tanh{\operatorname{th}}
\newcommand\ulie{\lie{u}}
\newcommand\upsi{{\underline{\psi}}}
\newcommand\A{\mathcal{A}}
\newcommand\Ad{\operatorname{Ad}}
\newcommand\Alt{\mathsf{Alt}}
\newcommand\Cset{\mathbb{C}}
\newcommand\Cycl{\mathop{\circlearrowleft}\limits}
\newcommand\DD{\mathbb{D}}
\newcommand\Dynl{\mathsfsl{Dynl}}
\newcommand\G{\mathcal{G}}
\newcommand\GG{\mathbb{G}}
\newcommand\Ker{\operatorname{\mathsf{Ker}}}
\newcommand\KK{\mathbb{K}}
\renewcommand\L{\mathcal{L}}
\newcommand\Lie{\operatorname{\mathsfsl{Lie}}}
\newcommand\M{\mathcal{M}}
\newcommand\Map{\mathsfsl{Map}}
\newcommand\Nset{\mathbb{N}}
\newcommand\Retire[1]{\relax}
\newcommand\Rset{\mathbb{R}}
\newcommand\REV{R^{\mathit{EV}}}
\renewcommand\S{\mathrm{S}}
\newcommand\Spectre{\operatorname{Sp}}
\newcommand\T{\operatorname{\mathsf{T}}}
\newcommand\Zset{\mathbb{Z}}
\newcommand\<{\mathopen\langle}
\renewcommand\>{\mathclose\rangle}

\begin{document}
\title{Quasi-bialgebras and dynamical $r$-matrices}
\author{Serge Parmentier}\email{serge@igd.univ-lyon1.fr}
\author{Romaric Pujol}\email{pujol@igd.univ-lyon1.fr}
\address{Institut Girard Desargues
  \newline
  \indent B\^at.~Jean Braconnier
  \newline
  \indent Universit\'e Claude Bernard -- Lyon I
  \newline
  \indent 21, avenue Claude Bernard
  \newline
  \indent F--69622 Villeurbanne Cedex.}
\date{February 25, 2004.}

\begin{abstract}
  We study the relationship between general dynamical Poisson groupoids
  and Lie quasi-bialgebras. For a class of Lie quasi-bialgebras $\G$
  naturally compatible with a reductive decomposition, we extend the
  description of the moduli space of classical dynamical $r$-matrices of
  Etingof and Schiffmann. We construct, in each gauge orbit, an explicit
  analytic representative $\lcan$. We translate the notion of duality
  for dynamical Poisson groupoids into a duality for Lie
  quasi-bialgebras. It is shown that duality maps the dynamical Poisson
  groupoid for $\lcan$ and $\G$ to the dynamical Poisson groupoid for
  $\lcan$ and the dual quasi-bialgebra $\G^\star$.
\end{abstract}
\maketitle

\section{Introduction}
The classical dynamical Yang--Baxter equation (CDYBE) for a pair
$(\g,\l\subset\g)$ of Lie algebras first appeared in~\cite{BDF,F}.
In~\cite{EV}, extending Drinfel\cprime d's classical work~\cite{D83},
this equation, supplemented by a condition of $\l$-equivariance, was
shown to coincide with the Jacobi identity for a natural Poisson bracket
on the trivial groupoid $U\times G\times U$. Here, $G$ is a Lie group
with $\Lie(G)=\g$ and $U\subset \l^*$ is an $\l$-invariant open set.

In~\cite{ES}, working on the formal disk $\DD$, Etingof and Schiffmann
have shown that the moduli space $\M(\DD,\l,\Omega)$ of formal solutions
of (CDYBE) for reductive triples $\g=\l\oplus\m$, $[\l,\m]\subset\m$,
where $\Omega\in\h\otimes\h\,\oplus\,\m\otimes\m$ is a symmetric
$\g$-invariant, is isomorphic to the algebraic variety
\begin{equation}
  \M_\Omega=
  \ensemble{\t\in(\wedge^2\m)^\l}
  {\<\t,\t\>+\<\Omega,\Omega\>\equiv0\mod\l}
\end{equation}
where $\<\,,\,\>$ is Drinfel\cprime d's bracket. Here, the moduli space
(first considered for $\Omega=0$ by Xu in~\cite{X}) is the orbit space
of solutions of (CDYBE) for the action of the group of equivariant maps
$\Map_0(\DD,G)^\l$ induced by formal base preserving, $\l$-equivariant
groupoid automorphisms of $\DD\times G\times\DD$. Their proof relies on
formal induction arguments and the equivariant Poincar\'e lemma together
with the existence of a canonical solution $\rAM$ of (CDYBE) on $\l$,
discovered by Alekseev and Meinrenken in~\cite{AM00}.

In~\cite{LP}, Poisson groupoid structures on $U\times G\times U$
compatible with the natural inclusion of the Hamiltonian unit $L\times
U$ (the so-called dynamical Poisson groupoids of~\cite{EV}) were
described. These brackets are given by pairs $(l,\varpi)$ where $l\colon
U\rightarrow \L(\g^*,\g)$ is a skew symmetric smooth map and $\varpi$ is
a $\g$-$1$-cocycle, and their Jacobi identity turns out to be equivalent
to
\begin{itemize}
\item There exists a $\varphi\in(\wedge^3\g)$ such that for all
  $\xi,\eta,\zeta\in \g^*$,
  \begin{equation}\tag{A}\label{eq:tagA}
    \<\xi\otimes\eta\otimes\zeta,\ad_x^{(3)}\varphi\>=
    \Cycl_{(\xi,\eta,\zeta)}
    \<\xi,\varpi_{\varpi_x\eta}\zeta\>
  \end{equation}
  and
  \begin{equation}\tag{B}\label{eq:tagB}
    \Cycl_{(\xi,\eta,\zeta)}\Big(\<\zeta,\d_pl(i^*\xi)\eta\>
    -\<\zeta,[l_p\xi,l_p\eta]\>
    -\<\zeta,\varpi_{l_p\xi}\eta\>\Big)=
    \<\xi\otimes\eta\otimes\zeta,\varphi\>
  \end{equation}
  for all $\xi,\,\eta,\,\zeta\in\g^*$, $x\in\g$, $p\in U$.
\item together with the $\l$-equivariance
  \begin{equation}\tag{C}\label{eq:tagC}
    \d_pl(\ad^*_zp)+\varpi_{iz}+\ad_{iz}l_p+l_p\ad^*_{iz}=0,
    \quad\forall z\in\l,\ p\in U
  \end{equation}
\end{itemize}
which, naturally, reduce to (CDYBE) for vanishing cocycle $\varpi$ and
$\g$-invariant $\varphi=\<\Omega,\Omega\>$.

The main purpose of the present work is to relate the above structures
to Lie quasi-bialgebras. Note that, for the case $\l=\g$, $\varpi=0$,
and $\varphi\in(\wedge^3\g)^\g$, such a link was also observed
in~\cite{EE}. First of all, we extend, with natural assumptions, the
Etingof--Schiffmann description of the moduli space in~\cite{ES} to the
equations~\eqref{eq:tagA}, \eqref{eq:tagB} and~\eqref{eq:tagC} with, as
principal new result, the construction, in terms of the associated Lie
quasi-bialgebra data, of a \emph{canonical analytic solution $\lcan$}
providing an explicit representative in each formal gauge class.
Secondly, we translate the notion of duality for Poisson
groupoids~\cite{W,MX2} in terms of a more algebraic \emph{duality for
  Lie quasi-bialgebras}, which up to (as yet formal) Poisson groupoid
automorphisms provides an explicit description of the Poisson groupoid
dual of the dynamical groupoid $U\times G\times U$.

Our analysis begins with the observation (see
Proposition~\ref{pr:Dynl->q-big}) that conditions~\eqref{eq:tagA},
\eqref{eq:tagB} and~\eqref{eq:tagC} above imply that the quadruple
$\G=(\g,[\,,\,],\varpi,\varphi)$ is a Lie quasi-bialgebra with
$\varpi_{\mid \l}$ exact. We then proceed to describe the moduli space
$\M$ of solutions of~\eqref{eq:tagB} and~\eqref{eq:tagC} for fixed
$\varpi$, $\varphi$ satisfying~\eqref{eq:tagA}. To begin with, in close
analogy with~\cite{ES}, we obtain (Corollary~\ref{cor:embedding}) an
embedding of $\M$ into the algebraic variety
\begin{equation*}
  \M_{\G,\l,\m}=\ensemble{\t\in(\wedge^2\m)^\l}{\varphi^\t\equiv0\mod\l},
\end{equation*}
where $\varphi^\t$ is the Drinfel'd twist of the associator
$\varphi\in\wedge^3\g$. Secondly for a class of Lie quasi-bialgebras
$\G$ canonically compatible (see Definition~\ref{def:Lie q-bialg
  compatible}) with a reductive decomposition $\g=\l\oplus\m$, we
construct the analytic solution $\lcan$ (see theorem~\ref{th:duality})
in terms of a generalization of $\rAM$ on $\l$ (also considered
in~\cite{EE}), and the \emph{adjoint action} of the double $\dlie$ of
the underlying Lie quasi-bialgebra. The use of $\lcan$ together with
Drinfel\cprime d twists, then shows (see
Corollary~\ref{cor:isomorphism}) that, for Lie quasi-bialgebras
compatible with the reductive decomposition $\g=\l\oplus\m$ the
embedding above is a bijection, providing, in particular, the explicit
analytic representative in each formal gauge orbit. Note that our
compatibility hypothesis on Lie quasi-bialgebras includes the class
considered in~\cite{ES}, for which the canonical solution $\lcan$
coincides, up to a twist, with the formal representative they
constructed.

For a contractible base $U$ (and, for large classes of examples, up to
covering, for arbitrary $U$ as well), it was shown in~\cite{LP} that
Poisson groupoid duality preserves the class of pairs $(U\times G\times
U,I)$ where $I\colon L\times U\rightarrow U\times G\times U$ is a
morphism of the Hamiltonian unit. However, the explicit expression of
the Poisson bracket of the dual pair $(U\times G'\times U,I')$ relies on
the knowledge of a non-canonical isomorphism (a so-called
trivialization) of the algebroid dual $A(U\times G\times U)^*\simeq
A(U\times G'\times U)$.

Our approach to duality begins with the construction of an explicit
trivialization of the algebroid dual for the canonical solution $\lcan$
(see Propositions~\ref{pr:LAB iso}, \ref{pr:flat connection & psi}, and
Theorem~\ref{th:trivialization}). It turns out that such an isomorphism
may be expressed solely in terms of the Drinfel\cprime d isomorphism
relating the doubles of the twisted pairs of Lie quasi-bialgebras $\G$
and $\G^{\lcan_q}$ together with the adjoint action of the double
$\dlie$ of $\G$.

Duality for Lie quasi-bialgebras is then defined (see
Definition~\ref{df:q-big dual}) as follows: let $\g=\l\oplus\m$ be a
reductive decomposition. If $\G=(\g,[\,,\,]\,\varpi,\varphi)$ is a Lie
quasi-bialgebra such that $\varpi_\l=0$ and $\varphi\equiv0\mod\l$, then
the dual $\G^\star$ is (up to relative signs) the Lie quasi-bialgebra
associated with the Manin quasi-triple
$(\dlie,\l\oplus\l^\perp,\m\oplus\m^\perp)$.

Our main duality assertion (see Theorem~\ref{th:duality}) then states
that the dual Poisson groupoid of the dynamical Poisson groupoid
associated with $\lcan$ for $\G$ is (isomorphic to) the
source-connected, simply-connected covering of the dynamical Poisson
groupoid associated with $\lcan$ for $\G^\star$. Note that this is
tantamount to saying that (up to covering) the dual of any dynamical
Poisson groupoid is dynamical if and only if the vertex algebra $\g$
admits a reductive decomposition $\g=\l\oplus\m$. This however will be
postponed to another publication~\cite{P}.

The paper is organized as follows. In section~\ref{sec:Liebialgebra}, we
recall some basic facts about Lie quasi-bialgebras with some additional
material needed for the rest of the paper. In section~\ref{sec:DPG and
  LQB} , we establish the relationship between dynamical Poisson
groupoids and Lie quasi-bialgebras. In section~\ref{sec:Gauge}, we adapt
the analysis of the moduli space in~\cite{ES} to the study of solutions
of~\eqref{eq:tagA}, \eqref{eq:tagB} and~\eqref{eq:tagC}. The brief
section~\ref{sec:cocom LQB} provides a formulation of the dynamical
$r$-matrix $\rAM$ in terms of the Lie quasi-bialgebra
$\G=(\g,[\,,\,],0,\varphi)$. In section~\ref{sec:com l-matrix}, we
construct the analytic representative $\lcan$ for Lie quasi-bialgebras
canonically compatible with a reductive decomposition. The last
section~\ref{sec:duality} is devoted to duality statements together with
some examples and a brief discussion on the link with duality of
symmetric spaces. We have collected some technical lemmas and proofs in
the Appendices.

\bigskip
\noindent\textbf{Acknowledgement.}
S.~Parmentier would like to thank A.~Alekseev for his invitation to the
University of Geneva at the early stages of this work. The financial
support of the Swiss National Science Foundation is gratefully
acknowledged.

\section{Lie quasi-bialgebras}
\label{sec:Liebialgebra}
In this section, we recall some basic facts about Lie quasi-bialgebras
(see~\cite{D90}, see also~\cite{AKS}).

\subsection{Notations}
When $E$ and $F$ are finite dimensional vector spaces over $\KK=\Rset$
or $\Cset$, we use the following notations:
\begin{itemize}
\item $E^*$ for the dual of $E$, and $\<\,,\,\>$ for the canonical
  pairing between $E$ and $E^*$,
\item $\L(E,F)$ for the set of linear maps from $E$ to $F$,
\item $f^*\in\L(F^*,E^*)$ for the adjoint of $f\in\L(E,F)$,
\item $\A(E^*,E)$ for the set of skew-symmetric linear maps from $E^*$
  to $E$.
\end{itemize}

Let $\g$ be a Lie algebra and let $G$ be a Lie group with $\Lie(G)=\g$.
In the sequel, Lie algebra and Lie group cocycles will always take value
in $\A(\g^*,\g)$ equipped with the adjoint action. Thus, a linear map
$\varpi\colon\g\to\A(\g^*,\g)$ is a Lie algebra $1$-cocycle if it
satisfies the following identity:
\begin{equation}
  \varpi_{[x,y]}=
  \ad_x\varpi_y+\varpi_y\ad^*_x-\ad_y\varpi_x-\varpi_x\ad^*_y
\end{equation}
for all $x,\,y\in\g$, and exact $1$-cocycles read as
$\varpi_x=\ad_x\t+\t\ad_x^*$, for some $\t\in\A(\g^*,\g)$. While a
smooth map $\pi\colon G\to\A(\g^*,\g)$ is a Lie group $1$-cocycle if it
satisfies the following identity:
\begin{equation}
  \pi_{gh}=\pi_g+\Ad_g\pi_h\Ad^*_g
\end{equation}
for all $g,\,h\in G$, and exact $1$-cocycles read as
$\pi_g=\Ad_g\t\Ad^*_g-\t$, where $\t\in\A(\g^*,\g)$.

Recall that $\varpi$ defined by $\varpi=\T_1\pi$ is a Lie algebra
$1$-cocycle. Moreover, if $G$ is connected and simply connected, then
Van Est's theorem (see \eg~\cite{G}) ensures that any Lie algebra
$1$-cocycle $\varpi$ may be uniquely lifted to a Lie group $1$-cocycle
such that $\varpi=\T_1\pi$.

The symbol ``$\Cycl_{(a_1,\ldots,a_n)}$'' means ``sum over cycling
permutations of $(a_1,\ldots,a_n)$''.

\subsection{Lie quasi-bialgebras}
\begin{defn}\label{df:doubledlie}
  Let $(\g,[\,,\,])$ be a Lie algebra, $\varpi\colon\g\to\A(\g^*,\g)$ a
  Lie algebra $1$-cocycle and $\varphi\in\wedge^3\g$. We say that the
  quadruple $\G=(\g,[\,,\,],\varpi,\varphi)$ is a \emph{Lie
    quasi-bialgebra} if $\dlie=\g\oplus\g^*$ together with the bracket
  $[\,,\,]_\dlie$
\begin{align}
  [x,y]_\dlie&=[x,y]\\*
  [x,\xi]_\dlie&=\varpi_x\xi-\ad^*_x\xi\\*
  [\xi,\eta]_\dlie&=\<\eta,\varpi_\bullet\xi\>+
  \<\xi\otimes\eta\otimes1,\varphi\>
\end{align}
for $x,\,y\in\g$ and $\xi,\,\eta\in\g^*$, is a Lie algebra. When
$\varpi=0$, the Lie quasi-bialgebra $(\g,[\,,\,],0,\varphi)$ is said to
be \emph{cocommutative}. The Lie algebra $(\dlie,[\,,\,]_\dlie)$ is
called the \emph{canonical double of the Lie quasi-bialgebra}
$(\g,[\,,\,],\varpi,\varphi)$.
\end{defn}

\begin{rem}
  Note that our sign convention for the associator $\varphi$ differs
  from that of Drinfel\cprime d in~\cite{D90}.
\end{rem}

The double $\dlie$ comes equipped with a non-degenerate invariant
symmetric bilinear form:
\begin{equation}
  (x+\xi,y+\eta)_\dlie=\<\xi,y\>+\<\eta,x\>
\end{equation}
for $x,\,y\in\g$ and $\xi,\,\eta\in\g^*$, for which $(\g,[\,,\,])$ is a
lagrangian (that is maximal isotropic) subalgebra of
$(\dlie,[\,,\,]_\dlie)$.

In practice, we will need the following:
\begin{prop}\label{prop:q-bialgebra}
  Let $(\g,[\,,\,])$ be a Lie algebra, let
  $\varpi\colon\g\to\A(\g^*,\g)$ be a Lie algebra $1$-cocycle, and
  $\varphi\in\wedge^3\g$. The quadruple $(\g,[\,,\,],\varpi,\varphi)$ is
  a Lie quasi-bialgebra if and only if the following two equations hold:
  \begin{align}
    \label{eq:Jacobi_1}
    \<\xi\otimes\eta\otimes\zeta,\ad^{(3)}_x\varphi\>-
    \Cycl_{(\xi,\eta,\zeta)}\<\xi,\varpi_{\varpi_x\eta}\zeta\>&=0\\*
    \Cycl_{(\xi,\eta,\zeta)}\Big(
    \<\<\eta,\varpi_\bullet\xi\>\otimes\zeta\otimes\theta,\varphi\>
    +\<\xi\otimes\eta\otimes\<\theta,\varpi_\bullet\zeta\>,\varphi\>\Big)&=0
  \end{align}
  for all $\xi,\,\eta,\,\zeta,\,\theta\in\g^*$ and $x\in\g$. In
  particular, the quadruple $(\g,[\,,\,],0,\varphi)$ is a Lie
  quasi-bialgebra if and only if $\varphi$ lies in
  $\left(\wedge^3\g\right)^\g$.
\end{prop}

\subsection{Manin pairs, Manin quasi-triples}
\begin{defn}
  Let $(\dlie,[\,,\,])$ be a Lie algebra together with a non-degenerate
  invariant symmetric bilinear form $(\,,\,)_\dlie$. We say that a pair
  $(\dlie,\g)$ is a \emph{Manin pair} if $\g$ is a lagrangian subalgebra
  of $\dlie$.

  We say that a triple $(\dlie,\g,\h)$ is a \emph{Manin quasi-triple} if
  the pair $(\dlie,\g)$ is a Manin pair, and if $\h$ is an isotropic
  complement of $\g$ in $\dlie$.
\end{defn}

Hence, if $(\g,[\,,\,],\varpi,\varphi)$ is a Lie quasi-bialgebra with
canonical double $\dlie$, the double $(\dlie,\g)$ is a Manin pair, and
the triple $(\dlie,\g,\g^*)$ is a Manin quasi-triple.

Conversely, let $(\dlie,\g,\h)$ be a Manin quasi-triple. Identifying
$\h$ with $\g^*$ by means of $(\, , \,)_\dlie$ provides a Lie
quasi-bialgebra structure on $\g$ denoted by $\G_{(\dlie,\g,\h)}$.  Its
cocycle $\varpi$ and associator $\varphi$ are explicitely given by
\begin{align}\label{eq:varpih}
  \varpi_x\xi&=\p_\g[x,\Omega^{-1}\xi]_\dlie\\*
  \label{eq:varphih}
  \<\xi\otimes\eta\otimes\zeta,\varphi\>&=
  \bigl(\Omega^{-1}\xi,[\Omega^{-1}\eta,
  \Omega^{-1}\zeta]_\dlie\bigr)_\dlie
\end{align}
where $\Omega$ is the identification $\Omega\colon\h\to\g^*$ given by
$(\,,\,)_\dlie)$.

\subsection{Twists}
Let $\G=(\g,[\,,\,],\varpi,\varphi)$ be a Lie quasi-bialgebra, with
canonical double $\dlie$. For an isotropic complement $\h$ of $\g$ in
$\dlie$, there exists a skew-symmetric linear map $\t\colon\g^*\to\g$
such that
\begin{equation*}
  \h=\ensemble{\t\xi+\xi}{\xi\in\g^*}
\end{equation*}
 The Lie
quasi-bialgebra induced by the Manin quasi-triple $(\dlie,\g,\h)$ is the
quadruple $(\g,[\,,\,],\varpi^\t,\varphi^\t)$ where:
\begin{align}\label{eq:varpi_twist}
  \varpi^\t_x&=\varpi_x+\ad_x\t+\t\ad_x^*\\*
  \label{eq:varphi_twist}
  \<\xi\otimes\eta\otimes\zeta,\varphi^\t\>&=
  \<\xi\otimes\eta\otimes\zeta,\varphi\>+
  \Cycl_{(\xi,\eta,\zeta)}
  \<\zeta,[\t\xi,\t\eta]+\varpi_{\t\xi}\eta\>
\end{align}
for $x\in\g$, $\xi,\,\eta,\,\zeta\in\g^*$. The Lie quasi-bialgebra
$(\g,[\,,\,],\varpi^\t,\varphi^\t)$ is called \emph{the twist of the Lie
  quasi-bialgebra $\G$ via $\t$}, and is denoted by $\G^\t$.  Note that
$\G$ is a Lie quasi-bialgebra if and only if $\G^\t$ is for any twist
$\t\in\A(\g^*,\g)$. Let $\dlie^\t$ be the double of $\g^\t$. The
following isomorphism of Drinfel\cprime d~\cite{D90}
\begin{equation}
  \appli{\tau_\t}{\dlie^\t}{\dlie}{x+\xi}{x+\t\xi+\xi}
\end{equation}
will play a crucial role in the sequel. Note that $\tau_\t$ preserves
the bilinear forms of $\dlie$ and $\dlie^\t$.  The inverse
$\tau_\t^{-1}$ of $\tau_\t$ is given by $\tau_\t^{-1}=\tau_{-\t}$, and,
if $\t'\in\A(\g^*,\g)$, then $\left(\G^\t\right)^{\t'}=\G^{\t+\t'}$.

We will use the following observation:
\begin{prop}\label{prop:prem_cond_twist}
  Let $\g$ be a Lie algebra, $\varpi\colon\g\to\A(\g^*,\g)$ a Lie
  algebra $1$-cocycle, and $\varphi\in\wedge^3\g$. For
  $\t\in\A(\g^*,\g)$, let $\varpi^\t$ and $\varphi^\t$ be as in
  equations~\eqref{eq:varpi_twist} and~\eqref{eq:varphi_twist}. Then the
  following equation holds:
  \begin{equation}
    \<\xi\otimes\eta\otimes\zeta,\ad^{(3)}_x\varphi\>-
    \Cycl_{(\xi,\eta,\zeta)}\<\xi,\varpi_{\varpi_x\eta}\zeta\>=
    \<\xi\otimes\eta\otimes\zeta,\ad^{(3)}_x\varphi^\t\>-
    \Cycl_{(\xi,\eta,\zeta)}\<\xi,\varpi^\t_{\varpi^\t_x\eta}\zeta\>
  \end{equation}
  for all $x\in\g$ and $\xi,\,\eta,\,\zeta\in\g^*$.
\end{prop}

\subsection{Lie quasi-bialgebra morphisms}
Let $\G_j=(\g_j,[\,,\,]_j,\varpi^j,\varphi^j)$, $j=1,\,2$ be two Lie
quasi-bialgebras, and $\upsi\colon\g_1\to\g_2$ a Lie algebra morphism.
We say that \emph{$\upsi$ is a Lie quasi-bialgebra morphism from $\G_1$
  to $\G_2$} if the following two conditions hold:
\begin{align}
  \upsi\varpi^1_x\upsi^*&=\varpi^2_{\upsi x}\quad\forall x\in\g_1\\*
  \upsi^{(3)}\varphi^1&=\varphi^2
\end{align}

The effect of twisting $\G_1$ via some $\t\in\A(\g_1^*,\g_1)$ is given
in the following proposition:
\begin{prop}
  Let $\upsi\colon\g_1\to\g_2$ be a Lie algebra morphism and let
  $\t\in\A(\g_1^*,\g_1)$. Set $\t'=\upsi\t\upsi^*$. Then the morphism
  $\upsi$ is a Lie quasi-bialgebra morphism from $\G_1$ to $\G_2$ if and
  only if $\upsi$ is a Lie quasi-bialgebra morphism from $\G_1^\t$ to
  $\G_2^{\t'}$.
\end{prop}

We also have the following lemma (for a proof, see
appendix~\ref{ap:pflem}):
\begin{lem}\label{lm:ad1 ad2}
  Let $\G_j=(\g_j,[\,,\,]_j,\varpi^j,\varphi^j)$, $j=1,\,2$ be two Lie
  quasi-bialgebras with double $\dlie^j$, and $\upsi$ a Lie
  quasi-bialgebra morphism from $\G_1$ to $\G_2$. Then the relations
  \begin{align}
    \label{eq:lm ad1 ad2 -- 1}
    \upsi\p_{\g_1}\left(\ad^1_{\upsi^*\xi}\right)^nu&=
    \p_{\g_2}\Bigl(\ad_\xi^2\Bigr)^n\upsi u\\*
    \label{eq:lm ad1 ad2 -- 1d}
    \p_{\g_1^*}\left(\ad^1_{\upsi^*\xi}\right)^n\upsi^*\eta&=
    \upsi^*\p_{\g_2^*}\Bigl(\ad_\xi^2\Bigr)^n\eta\\*
    \label{eq:lm ad1 ad2 -- 2}
    \p_{\g_1^*}\left(\ad^1_{\upsi^*\xi}\right)^nu&=
    \upsi^*\p_{\g_2^*}\Bigl(\ad_\xi^2\Bigr)^n\upsi u\\*
    \label{eq:lm ad1 ad2 -- 3}
    \upsi\p_{\g_1}\left(\ad^1_{\upsi^*\xi}\right)^n\upsi^*\eta&=
    \p_{\g_2}\Bigl(\ad_\xi^2\Bigr)^n\eta
  \end{align}
  hold for all $n\in\Nset$ and for all $u\in\g_1$,
  $\xi,\,\eta\in\g_2^*$.
\end{lem}

\subsection{Lie quasi-bialgebras obtained from one another}
\label{sec:LQBOFA}
Let $\G=(\g,[\,,\,],\varpi,\varphi)$ be a Lie quasi-bialgebra. Let
$\upsi$ be an automorphism of the Lie algebra $(\g,[\,,\,])$, and set
\begin{align}
  \label{eq:varpipsi}
  \varpi^\upsi_x=\upsi\varpi_{\upsi^{-1}x}\upsi^*\\*
  \label{eq:varphipsi}
  \varphi^\upsi=\upsi^{(3)}\varphi
\end{align}
Then $\G^\upsi=(\g,[\,,\,],\varpi^\upsi,\varphi^\upsi)$ is a Lie
quasi-bialgebra, and $\upsi$ is a Lie quasi-bialgebra isomorphism from
$\G$ to $\G^\upsi$. More generally, let $w$ be an automorphism of the
vector space $\g$, equip $\g$ with the bracket
\begin{equation}
  [x,y]^w=w^{-1}[wx,wy].
\end{equation}
and set $\varphi^w=w^{(3)}\varphi$ and $\varpi_x^w=w\varpi_{w^{-1}x}w^*$
for $x\in\g$. Then, the quadruple
$\G^w=(\g,[\,,\,]^w,\varpi^w,\varphi^w)$ is a Lie quasi-bialgebra such
that $w$ is a Lie quasi-bialgebra isomorphism between $\G$ and $\G^w$.

The Lie quasi-bialgebra $\G^-=(\g,[\,,\,],-\varpi,\varphi)$ is called
the \emph{inversion of the Lie quasi-bialgebra $\G$}. Obviously,
$(\G^-)^-=\G$. If we denote by $\dlie$ and $\dlie^-$ the double of $\G$
and $\G^-$ respectively, then the map $J\colon\dlie\to\dlie^-$ defined
by $J(x+\xi)=x-\xi$ for all $x\in\g$ and $\xi\in\g^*$ is a Lie algebra
isomorphism.

\subsection{The adjoint action}
Let $\G=(\g,[\,,\,],\varpi,\varphi)$ be a Lie quasi-bialgebra with
canonical double $\dlie$, let $D$ be the connected, simply-connected Lie
group with Lie algebra $\dlie$, let $G$ be the connected Lie subgroup of
$D$ with Lie algebra $\g$, and let $\pi\colon G\to\A(\g^*,\g)$ be the
Lie group $1$-cocycle integrating the Lie algebra $1$-cocycle $\varpi$.
Denote by $\Ad^D$ the adjoint action of $D$ on its Lie algebra $\dlie$.
For any $x\in\g$, $\xi\in\g^*$ and $g\in G$, one has:
\begin{equation}\label{eq:AdD}
  \Ad^D_g(x+\xi)=\Ad_gx+\pi_g\Ad^*_{g^{-1}}\xi+\Ad^*_{g^{-1}}\xi
\end{equation}
Indeed, it is easy to show that $\Ad^D_g\Ad^D_{g'}=\Ad^D_{gg'}$ for all
$g,\,g'\in G$, and that $\left.\frac\d{\d
    t}\right\vert_{t=0}\Ad^D_{\e^{tu}}(x+\xi)=\ad^\dlie_u(x+\xi)$ for
all $u\in\g$ and $x\in\g$, $\xi\in\g^*$.

\section{Dynamical Poisson groupoids and Lie quasi-bialgebras}
\label{sec:DPG and LQB}
For further informations on dynamical Poisson groupoids, see~\cite{EV}
and~\cite{LP}.

\subsection{Lie quasi-bialgebra associated with a trivial Poisson
  groupoid}
Let $G$ be a connected Lie group with Lie algebra $\g$. For any point
$x\in G$, we denote by $D_xf\in\g^*$ and $D'_xf\in\g^*$ the right and
left derivatives at $x$:
\begin{align}
  D_xf(u)&=\left.\frac\d{\d t}\right\vert_{t=0}f(\e^{tu}x)\\*
  \label{eq:D'_xf}
  D'_xf(u)&=\left.\frac\d{\d t}\right\vert_{t=0}f(x\e^{tu})
\end{align}
for all $u\in\g$. Let $L$ be a connected Lie subgroup of $G$ with Lie
algebra $\l$, and $U$ an $\Ad^*_L$-invariant open subset in $\l^*$. We
will denote the inclusion by $i\colon\l\to\g$. Consider the trivial Lie
groupoid $\GG=U\times G\times U$ with multiplication:
\begin{equation}
  (p,x,q)(q,y,r)=(p,xy,r)
\end{equation}
We say that a multiplicative Poisson bracket on $\GG$ is
\emph{dynamical} if it is of the form:
\begin{equation}\label{eq:Dynamical_bracket}
\begin{aligned}
  \{f,g\}_{(p,x,q)}=&\,\<p,[\delta f,\delta g]_\l\>
  -\<q,[\delta'f,\delta'g]_\l\>\\*
  &\quad-\<Dg,i\delta f\>-\<D'g,i\delta'f\>\\*
  &\quad+\<Df,i\delta g\>+\<D'f,i\delta'g\>\\*
  &\quad-\<Df,l_pDg\>+\<Df,\pi_xDg\>+\<D'f,l_qD'g\>
\end{aligned}
\end{equation}
where $l\colon U\to\A(\g^*,\g)$ is a smooth map, and $\pi\colon
G\to\A(\g^*,\g)$ is a group $1$-cocycle. In this equation, $\delta f$
and $\delta'f$ denote the derivatives of $f$ with respect to the first
and second $U$ factors, $Df$ and $D'f$ denote the right and left
derivatives of $f$ with respect to the $G$ factor, and all derivatives
are evaluated at $(p,x,q)$. Denote by $\varpi=\T_1\pi$ the Lie algebra
$1$-cocycle associated with $\pi$.

The map $P$ from $\GG$ to $\A(\g^*,\g)$ defined by:
\begin{equation}
  P_{(p,x,q)}=-l_p+\pi_x+\Ad_xl_q\Ad^*_x
\end{equation}
is called the \emph{groupoid cocycle associated with the dynamical Poisson
  bracket~\eqref{eq:Dynamical_bracket}}.

Using theorem~2.2.5.~of~\cite{LP}, it may be shown that the Jacobi
identity for a bracket of this type is equivalent to the following two
conditions:
\begin{itemize}
\item There exists a $\varphi\in(\wedge^3\g)$ such that for all
  $\xi,\,\eta,\,\zeta\in\g^*$:
  \begin{equation}\label{eq:phi-eqg}
    \<\xi\otimes\eta\otimes\zeta,\ad_x^{(3)}\varphi\>=
    \Cycl_{(\xi,\eta,\zeta)}
    \<\xi,\varpi_{\varpi_x\eta}\zeta\>
  \end{equation}
  and for all $p\in U$ and $\xi,\,\eta,\,\zeta\in\g^*$:
  \begin{equation}\label{eq:CDYBEg}
    \Cycl_{(\xi,\eta,\zeta)}\Big(\<\zeta,\d_pl(i^*\xi)\eta\>
    -\<\zeta,[l_p\xi,l_p\eta]\>
    -\<\zeta,\varpi_{l_p\xi}\eta\>\Big)=
    \<\xi\otimes\eta\otimes\zeta,\varphi\>
  \end{equation}
\item For all $p\in U$ and $z\in\l$:
  \begin{equation}\label{eq:l-eqg}
    \d_pl(\ad^*_zp)+\varpi_{iz}+\ad_{iz}l_p+l_p\ad^*_{iz}=0
  \end{equation}
\end{itemize}

Equation~\eqref{eq:CDYBEg} can be seen as a generalization of the
modified classical dynamical Yang--Baxter equation to which it reduces
when $\varpi=0$. Equation~\eqref{eq:phi-eqg} is exactly
equation~\eqref{eq:Jacobi_1}, and equation~\eqref{eq:l-eqg} is a
generalization of the $\l$-equivariance of the map $l$.

\begin{rem}
  Equation~\eqref{eq:CDYBEg} may also be written:
  \begin{multline}\label{eq:CDYBEg alt}
    \d_pl(i^*\xi)\eta-\d_pl(i^*\eta)\xi-i\d_p\<\xi,l_\bullet\eta\>
    -[l_p\xi,l_p\eta]-l_p\ad^*_{l_p\xi}\eta+l_p\ad^*_{l_p\eta}\xi\\*
    -\varpi_{l_p\xi}\eta+\varpi_{l_p\eta}\xi+
    \<\xi,\varpi_{l_p\bullet}\eta\>=\<\xi\otimes\eta\otimes1,\varphi\>
  \end{multline}
  for all $\xi,\,\eta\in\g^*$.
\end{rem}

Now, by a result from~\cite{LP}, we know that, for a contractible base
$U$, the dual (see section~\ref{sec:duality} below, for more explicit
information about duality) of the Poisson groupoid $\GG$ with Poisson
bracket~\eqref{eq:Dynamical_bracket} is still a trivial Poisson groupoid
(not necessarily dynamical, though). Its vertex Lie group
$G^\star_{q_0}$ is the connected, simply connected Lie group with Lie
algebra (isomorphic to) the vector space
\begin{equation}\label{eq:dual_of_g}
  \g_{q_0}^\star=\ensemble{i(z)+\xi\in i(\l)\oplus\g^*}
  {i^*\xi=\ad^*_zq_0}\subset\g\oplus\g^*
\end{equation}
 for some $q_0\in U$, together with the Lie bracket:
\begin{equation}\label{eq:bracket_on_gstar}
  \begin{aligned}\relax
    [i(z)+\xi,i(z')+\xi']^\star_{q_0}&=\bigl(i([z,z'])+\varpi_{i(z)}\xi'
    +\ad_{i(z)}l_{q_0}\xi' +l_{q_0}\ad_{i(z)}^*\xi'\\*
    &\quad-\varpi_{i(z')}\xi-\ad_{i(z')}l_{q_0}
    \xi-l_{q_0}\ad_{i(z')}^*\xi\\*
    &\quad+[l_{q_0}\xi,l_{q_0}\xi']+l_{q_0}\ad^*_{l_{q_0}\xi}\xi'
    -l_{q_0}\ad^*_{l_{q_0}\xi'}\xi\\*
    &\quad+\varpi_{l_{q_0}\xi}\xi'-\varpi_{l_{q_0}\xi'}\xi
    -\<\xi,\varpi_{l_{q_0}\bullet}\xi'\>
    +\<\xi\otimes \xi'\otimes1,\varphi\>,\\*
    &\quad-\ad^*_{i(z)}\xi'+\ad^*_{i(z')}\xi-\<\xi,\varpi_\bullet\xi'\>
    -\ad^*_{l_{q_0}\xi}\xi'+\ad^*_{l_{q_0}\xi'}\xi\bigr)
  \end{aligned}
\end{equation}
for all $i(z)+\xi,\,i(z')+\xi'\in\g^\star_{q_0}$. Note that the Lie
algebras $(\g_{q_0}^\star,[\,,\,]^\star_{q_0})$ are all isomorphic when
$q_0$ ranges over $U$.

We start with a lemma which relates a solution $l$ of
equation~\eqref{eq:CDYBEg} to its translation $l'=l-\t$ by an element
$-\t\in\A(\g^*,\g)$:
\begin{lem}\label{lm:CDYBEg twist}
  Let $\t\in\A(\g^*,\g)$. Set $l'_p=l_p-\t$ for any $p\in U$. Then
  equation~\eqref{eq:CDYBEg} is satisfied for all
  $\xi,\,\eta,\,\zeta\in\g^*$ and $p\in U$ if and only if the following
  equation is statisfied:
  \begin{equation}\label{eq:CDYBEg_twisted}
    \Cycl_{(\xi,\eta,\zeta)}\Big(\<\zeta,\d_pl'(i^*\xi)\eta\>
    -\<\zeta,[l'_p\xi,l'_p\eta]\>
    -\<\zeta,\varpi^\t_{l'_p\xi}\eta\>\Big)=
    \<\xi\otimes\eta\otimes\eta,\varphi^\t\>
  \end{equation}
  for all $\xi,\,\eta,\,\zeta\in\g^*$ and $p\in U$, where $\varpi^\t$
  and $\varphi^\t$ are defined by equations~\eqref{eq:varpi_twist}
  and~\eqref{eq:varphi_twist} (even though we don't know yet that
  $(\g,[\,,\,],\varpi,\varphi)$ is a Lie quasi-bialgebra).
\end{lem}
\begin{proof}
  Straightforward computation using equations~\eqref{eq:varpi_twist}
  and~\eqref{eq:varphi_twist}.
\end{proof}

As an immediate consequence we can write equation~\eqref{eq:CDYBEg} as
\begin{equation}\label{eq:varphi_twisted}
  \<\xi\otimes\eta\otimes\zeta,\varphi^{l_p}\>=
  \Cycl_{(\xi,\eta,\zeta)}\<\zeta,\d_pl(i^*\xi)\eta\>
\end{equation}
for all $p\in U$ and $\xi,\,\eta,\,\zeta\in\g^*$. Thus, if $l$ satisfies
equation~\eqref{eq:CDYBEg} on $U$, then $\varphi^{l_p}\equiv0\mod\l$,
$\forall p\in U$.

We now come to a proposition which is basic for our subsequent analysis.
\begin{prop}\label{pr:Dynl->q-big}
  Let $q_0$ be a point in $U$. If equations~\eqref{eq:phi-eqg},
\eqref{eq:CDYBEg},
 and~\eqref{eq:l-eqg} are satisfied, then the
  quadruple $\G^{q_0}=(\g,[\,,\,],\varpi^{l_{q_0}},\varphi^{l_{q_0}})$
  is a Lie quasi-bialgebra. Moreover, the Lie algebra $\g_{q_0}^\star$
  defined by equations~\eqref{eq:dual_of_g}
  and~\eqref{eq:bracket_on_gstar} is a Lagrangian subalgebra of the
  canonical double $\dlie^{q_0}$ of $\G^{q_0}$.
\end{prop}
Before proving proposition~\ref{pr:Dynl->q-big}, we state the following
auxiliary result:
\begin{lem}\label{lem:gnbwsjnklgkl}
  Let $l$ be a solution of~\eqref{eq:CDYBEg} and~\eqref{eq:l-eqg} on $U$
  and set $l'=l-l_{q_0}$. Then for all
  $\xi,\,\eta,\,\zeta,\,\theta\in\g^*$, the two following equations
  hold:
  \begin{align}
    \label{eq:blabla1}
    \<\theta,
    \d_{q_0}l'\big(i^*\<\eta,\varpi^{l_{q_0}}_\bullet\xi\>\big)\zeta\>&=
    \<\xi,\d_{q_0}l'
    \big(i^*\<\theta,\varpi^{l_{q_0}}_\bullet\zeta\>\big)
    \eta\>\\*
    \label{eq:derivee_CDYBEg}
    \Cycl_{(\xi,\eta,\zeta)}
    \<\zeta,\d^2_{q_0}l'(i^*\xi,i^*\theta)\eta\>&=
    \Cycl_{(\xi,\eta,\zeta)}
    \<\zeta,\varpi^{l_{q_0}}_{\d_{q_0}l'(i^*\theta)\xi}\eta\>
  \end{align}
\end{lem}
\begin{proof}
  Equation~\eqref{eq:blabla1} is a consequence of
  equation~\eqref{eq:l-eqg}, and equation~\eqref{eq:derivee_CDYBEg} is
  the derivative of equation~\eqref{eq:CDYBEg_twisted} in the direction
  $i^*\theta$, evaluated at $q_0$.
\end{proof}

\begin{proof}[Proof of proposition~\ref{pr:Dynl->q-big}]
  According to proposition~\ref{prop:prem_cond_twist}, the first
  condition of proposition~\ref{prop:q-bialgebra} is satisfied, since
  equation~\eqref{eq:phi-eqg} holds, so it only remains to show the
  second condition of proposition~\ref{prop:q-bialgebra}.

  According to equality~\eqref{eq:varphi_twisted}, the second condition
  of proposition~\ref{prop:q-bialgebra} for
  $(\g,[\,,\,],\varpi^{l_{q_0}},\varphi^{l_{q_0}})$ reads:
  \begin{multline}
    \Cycl_{(\xi,\eta,\zeta)}\Big(
    \<\<\eta,\varpi^{l_{q_0}}_\bullet\xi\>\otimes\zeta\otimes\theta,
    \varphi^{l_{q_0}}\>
    +\<\xi\otimes\eta\otimes\<\theta,\varpi^{l_{q_0}}_\bullet\zeta\>,
    \varphi^{l_{q_0}}\>\Big)=\\*
    \Cycl_{(\xi,\eta,\zeta)}
    \<\eta,\varpi^{l_{q_0}}_{\d_{q_0}l'(i^*\zeta)\theta}\xi\>+
    \<\xi,\varpi^{l_{q_0}}_{\d_{q_0}l'(i^*\theta)\zeta}\eta\>+ \<\theta,
    \d_{q_0}l'\left(i^*\<\eta,
      \varpi^{l_{q_0}}_\bullet\xi\>\right)\zeta\>+\\*
    \qquad\<\theta,\varpi^{l_{q_0}}_{\d_{q_0}l'(i^*\xi)\eta}\zeta\>+
    \<\zeta,\varpi^{l_{q_0}}_{\d_{q_0}l'(i^*\eta)\xi}\theta\>+ \<\eta,
    \d_{q_0}l'
    \left(i^*\<\theta,\varpi^{l_{q_0}}_\bullet\zeta\>\right)\xi\>
  \end{multline}
  which vanishes by lemma~\ref{lem:gnbwsjnklgkl} and Schwarz' lemma.
  Hence, the quadruple
  $\G^{q_0}=(\g,[\,,\,],\varpi^{l_{q_0}},\varphi^{l_{q_0}})$ is a Lie
  quasi-bialgebra.

  It is clear from equation~\eqref{eq:bracket_on_gstar} that
  $\g_{q_0}^\star$ is a Lie subalgebra of the canonical double
  $\dlie^{q_0}$ of $\G^{q_0}$, and a simple verification shows that it
  is lagrangian.
\end{proof}

As mentioned above, beware that the dual Poisson groupoid of a dynamical
Poisson groupoid is not dynamical in general, and even if it is, the Lie
quasi-bialgebra on $\g^\star_0$, which the dual Poisson groupoid is
associated with, is not directly obtained from the Manin pair
$(\dlie,\g^\star_0)$.

We now introduce the following definitions:
\begin{defn}
  Let $\G=(\g,[\,,\,],\varpi,\varphi)$ be a Lie quasi-bialgebra, $\l$ a
  Lie subalgebra of $\g$, and $U\subset\l^*$ an $\Ad^*_L$-invariant open
  subset.
  \begin{enumerate}
  \item We say that a smooth map $l\colon U\to\A(\g^*,\g)$ is a
    \emph{dynamical $\ell$-matrix on $U$ associated with the Lie
      quasi-bialgebra $\G$} if it satisfies equations~\eqref{eq:CDYBEg}
    and~\eqref{eq:l-eqg}. A dynamical $\ell$-matrix associated with a
    cocommutative Lie quasi-bialgebra is called a \emph{dynamical
      $r$-matrix.}
  \item Let $q\in\l^*$, and let $\DD_q\subset\l^*$ be the formal
    neighborhood of $q$. We say that a (formal) map
    $l\colon\DD_q\to\A(\g^*,\g)$ is a \emph{formal dynamical
      $\ell$-matrix at $q$ associated with the Lie quasi-bialgebra $\G$}
    if it satisfies equations~\eqref{eq:CDYBEg} and~\eqref{eq:l-eqg}
    formally.
  \item We denote by $\Dynl(U,\G)$ the set of dynamical $\ell$-matrices
    on $U$ associated with the Lie quasi-bialgebra $\G$, and by
    $\Dynl(\DD_q,\G)$ the set of formal dynamical $\ell$-matrices at $q$
    associated with the Lie quasi-bialgebra $\G$.
  \end{enumerate}
\end{defn}

Lemma~\ref{lm:CDYBEg twist} has the following interpretation:
\begin{prop}\label{pr:Dynl twist}
  With these notations, $\forall\t\in\A(\g^*,\g)$,
  \begin{equation*}
    \begin{aligned}
      \Dynl(U,\G^\t)&=\Dynl(U,\G)-\t\\*
      \Dynl(\DD_q,\G^\t)&=\Dynl(\DD_q,\G)-\t
    \end{aligned}
  \end{equation*}
\end{prop}

>From now on, we assume that $0\in U$, we set $\DD=\DD_0$, and we define
\begin{equation*}
  \begin{aligned}
    \Dynl_0(U,\G)&=\ensemble{l\in\Dynl(U,\G)}{l_0=0}\\*
    \Dynl_0(\DD,\G)&=\ensemble{l\in\Dynl(\DD,\G)}{l_0=0}
  \end{aligned}
\end{equation*}
Using proposition~\ref{pr:Dynl twist}, we get:
\begin{equation*}
  \begin{aligned}
    \Dynl(U,\G)&=\bigcup_{\t\in\A(\g^*,\g)}\Dynl_0(U,\G^\t)+\t\\*
    \Dynl(\DD,\G)&=\bigcup_{\t\in\A(\g^*,\g)}\Dynl_0(\DD,\G^\t)+\t
  \end{aligned}
\end{equation*}

\begin{rem}\label{rk:non empty dynl0}
  For $\Dynl(U,\G)$ to be non empty, it is necessary that there exists a
  twist $\t\in\A(\g^*,\g)$ such that $\varpi^\t_\l=0$ and that
  $\varphi^\t\equiv0\mod\l$. Indeed, if $l\in\Dynl(U,\G)$, then
  $l'=l-l_0\in\Dynl_0(U,\G^{l_0})$, as shown by proposition~\ref{pr:Dynl
    twist}. Now lemma~\ref{lm:CDYBEg twist} implies that
  $\varphi^{l_0}\equiv0\mod\l$ and equation~\eqref{eq:l-eqg} implies
  that $\varpi^{l_0}_\l=0$. Obviously, the same holds for
  $\Dynl(\DD,\G)$.
\end{rem}

\begin{rem}
  Let $q_0\in\l^*$ such that $\ad^*_zq_0=0$ for all $z\in\l$, and let
  $l\in\Dynl(U,\G)$. Then the map $l'\colon U+q_0\to\A(\g^*,\g)$ defined
  by $l'_{p}=l_{p-q_0}$ lies in $\Dynl(U+q_0,\G)$.
\end{rem}

\section{Gauge transformations}
\label{sec:Gauge}
In this section we recall the action of the gauge group on dynamical
$\ell$-matrices which was introduced in~\cite{EV} for dynamical
$r$-matrices. Also, we describe the associated moduli space, following
the scheme of~\cite{ES}.

For any subset $A$ of the vector space $E$, we denote by $A^\perp$ the
orthogonal space to $A$:
\begin{equation*}
  A^\perp=\ensemble{v\in E^*}{\<v,a\>=0,\,\forall a\in A}.
\end{equation*}

We shall denote by $\ell_x$ and $r_x$ the left and right action of a Lie
group $G$ on its tangent bundle associated with the left and right
multiplications of $G$ on itself.

\subsection{Trivial groupoid morphisms}
Let $\GG_1=U\times G_1\times U$ and $\GG_2=U\times G_2\times U$ be two
trivial Lie groupoids over the same base $U$ which is assumed to contain
$0$. Let $\Psi\colon\GG_1\to\GG_2$ be a base preserving groupoid
morphism.
\begin{prop}\label{pr:groupoid morphisms}
  The morphism $\Psi$ has the form:
  \begin{equation}
    \Psi(p,x,q)=(p,\sigma_p\psi(x)\sigma_q^{-1},q)
  \end{equation}
  for all $(p,x,q)\in\GG_1$, where $\sigma\colon U\to G_2$ is a smooth
  map satisfying $\sigma_0=1$, and $\psi\colon G_1\to G_2$ is a Lie
  group morphism.
\end{prop}
\begin{proof}
  The most general form for a base preserving map
  $\Psi\colon\GG_1\to\GG_2$ is:
  \begin{equation}
    \Psi(p,x,q)=(p,\psi_{p,q}(x),q)
  \end{equation}
  where $\psi_{p,q}\colon G_1\to G_2$. We set $\psi=\psi_{0,0}$ and
  $\sigma_p=\psi_{p,0}(1)$. If $\Psi$ is a groupoid morphism, then
  $\psi_{p,q}(x)=\psi_{p,0}(1)\psi_{0,0}(x)\psi_{0,q}(1)$ and
  $\psi_{0,q}(x)=\psi_{q,0}(x^{-1})^{-1}$ for all $x\in G$. Thus, $\psi$
  is a Lie group morphism, $\psi_{p,q}(x)=\sigma_p\psi(x)\sigma_q^{-1}$
  and $\sigma_0=1$.
\end{proof}

Let $L$ be a connected Lie subgroup of both $G_1$ and $G_2$, with Lie
algebra $\l$. If $U$ is an $\Ad^*_L$-invariant subset of $\l^*$, there
are two actions of $L$ on $\GG_k$, $k=1,\,2$, namely:
\begin{itemize}
\item A left action: $h\cdot(p,x,q)=(\Ad^*_{h^{-1}}p,hx,q)$,
\item and a right action: $(p,x,q)\cdot h=(p,xh,\Ad^*_hq)$.
\end{itemize}
By definition, the groupoid morphism $\Psi$ is said to be
$L$-biequivariant if and only if it is equivariant for both left and
right actions of $L$, that is if and only if
\begin{equation}
  \sigma\left(\Ad^*_{h^{-1}}p\right)=h\sigma_p\psi(h)^{-1}
\end{equation}
for all $h\in L$ and $p\in U$. Since $L$ is connected and $0\in U$, this
condition is also equivalent to its infinitesimal version:
\begin{align}
  \label{eq:psi=id on l}
  \upsi z&=z\\*
  \label{eq:sigma l-eq}
  r_{\sigma_p^{-1}}(\T_p\sigma)\ad^*_zp&=
  \Ad_{\sigma_p}z-z
\end{align}
for all $z\in\l$ and $p\in U$, where $\upsi=\T_1\psi\colon\g_1\to\g_2$
is the Lie algebra morphism associated to the Lie group morphism $\psi$.

Now, Let $P^1$ and $P^2$ be groupoid cocycles on $\GG_1$ and $\GG_2$
associated with some dynamical Poisson brackets on $\GG_1$ and $\GG_2$.
For $f\in C^\infty(\GG_2)$ and $X=(p,x,q)$ a point in $\GG_1$, a
direct calculation yields:
\begin{align}
  \delta_X(f\circ\Psi)&=\delta_{\Psi(X)}f
  +(\T_p\sigma)^*r^*_{\sigma_p^{-1}}D_{\Psi(X)}f\\*
  D_X(f\circ\Psi)&=\upsi^*\Ad^*_{\sigma_p}D_{\Psi(X)}f\\*
  D'_X(f\circ\Psi)&=\upsi^*\Ad^*_{\sigma_q}D'_{\Psi(X)}f\\*
  \delta'_X(f\circ\Psi)&=\delta'_{\Psi(X)}f
  -(\T_q\sigma)^*r^*_{\sigma_q^{-1}}D'_{\Psi(X)}f
\end{align}
Using these equations and equations~\eqref{eq:psi=id on l}
and~\eqref{eq:sigma l-eq}, one can show the following:
\begin{prop}\label{pr:dyn Poisson groupoid morphism}
  The groupoid morphism $\Psi$ is a Poisson groupoid morphism if and
  only if the two following conditions hold:
  \begin{itemize}
  \item The groupoid morphism $\Psi$ is $L$-biequivariant,
  \item The equation
    \begin{equation}\label{eq:action on grpd cocyles}
      \Ad_{\sigma_p}\upsi P^1_X\upsi^*\Ad_{\sigma_p}^*
      +\Theta^\Psi_{\Psi(X)}=P^2_{\Psi(X)}
    \end{equation}
    is satisfied for all $X=(p,x,q)\in\GG_1$, where
    \begin{align}
      \Theta^\Psi_{(p,y,q)}&=
      \Ad_y\theta^\Psi_q\Ad^*_y-\theta^\Psi_p\\*
      \label{eq:theta}
      \theta^\Psi_p&=r_{\sigma_p^{-1}}(\T_p\sigma)i^*\Ad_{\sigma_p}^*
      -(\T_p\sigma)^*r^*_{\sigma_p^{-1}}
    \end{align}
    for all $(p,y,q)\in\GG_2$.
  \end{itemize}
\end{prop}

Notice that $\theta^\Psi$, as defined by equation~\eqref{eq:theta} is
skew-symmetric (use equation~\eqref{eq:sigma l-eq}), and that
$\Theta^\Psi$ is an exact groupoid $1$-cocycle.

For $j=1,\,2$, write
$P^j_{(p,x_j,q)}=-l^j_p+\pi^j_{x_j}+\Ad_{x_j}l^j_q\Ad^*_{x_j}$ for all
$(p,x_j,q)\in\GG_j$. Without loss of generality, we may assume that
$l^2_0=\theta^\Psi_0+\upsi l^1_0\upsi^*$ (this is done by translating
$l^2$, while adding an exact group-cocycle to $\pi^2$). Denote by
$\G_j=(\g_j,[\,,\,],\varpi^j,\varphi^j)$ the Lie quasi-bialgebras
associated with $l^j$ and the Poisson bracket on $\GG_j$.
\begin{prop}\label{pr:gpdmorph-->qbigmorph}
  With this convention, the following equation holds for all $p\in U$:
  \begin{equation}\label{eq:l^2^Psi}
    l^2_p=\Ad_{\sigma_p}\upsi l_p^1\upsi^*\Ad^*_{\sigma_p}
    +\theta^\Psi_p+\pi^2_{\sigma_p}
  \end{equation}
  and the Lie algebra morphism $\upsi$ is a Lie quasi-bialgebra morphism
  from $\G_1$ to $\G_2$.
\end{prop}
\begin{proof}
  Evaluating equation~\eqref{eq:action on grpd cocyles} at $X=(0,1,p)$
  for some $p\in U$ yields equation~\eqref{eq:l^2^Psi}, and evaluating
  equation~\eqref{eq:action on grpd cocyles} at $X=(0,x,0)$ for $x\in
  G_1$ yields:
  \begin{equation}
    \upsi\pi^1_x\upsi^*=\pi^2_{\psi x}
  \end{equation}
  There exists an open subset $U'$ of $U$, with $0\in U'$, and a smooth
  map $\Sigma\colon U'\to\g$ satisfying $\Sigma_0=0$ such that
  $\sigma_p=\e^{\Sigma_p}$ for all $p\in U'$. Set
  $A=\d_0\Sigma\in\L(\l^*,\g_2)$. Notice that $A^*$ takes values in $\l$
  and that $A$ is $\l$-equivariant, since $\sigma$ is. Using the
  differential of the exponential map (see section~\ref{sec:diff_exp} in
  the appendix), one computes (here and below, $i^*$ stands for
  $i_2^*$):
  \begin{equation}
    \theta^\Psi_p=
    \frac{\Ad_{\e^{\Sigma_p}}-1}{\ad_{\Sigma_p}}(\d_p\Sigma)i^*
    \Ad^*_{\e^{\Sigma_p}}
    -(\d_p\Sigma)^*\frac{\Ad^*_{\e^{\Sigma_p}}-1}{\ad^*_{\Sigma_p}}
  \end{equation}
  so that
  \begin{align}\label{eq:0^Psi}
    \theta^\Psi_0&=Ai^*-A^*\\*
    \label{eq:dtheta0}
    \begin{split}
      \<\eta,\d_0\theta^\Psi(\alpha)\xi\>&=
      \frac12\<\eta,[A\alpha,Ai^*\xi]\>-\frac12\<\xi,[A\alpha,Ai^*\eta]\>
      +\<\eta,Ai^*\ad^*_{A\alpha}\xi\>\\*
      &\qquad
      +\<\eta,\d_0^2\Sigma(\alpha,i^*\xi)\>
      -\<\xi,\d_0^2\Sigma(\alpha,i^*\eta)\>
    \end{split}
  \end{align}
  for all $\alpha\in\l^*$, and $\xi,\,\eta\in\g^*$. Now, by
  equation~\eqref{eq:CDYBEg}, we have for all
  $\xi,\,\eta,\,\zeta\in\g^*_2$:
  \begin{equation}\label{eq:CDYBEgPsi}
    \<\xi\otimes\eta\otimes\zeta,\varphi^2\>=
    \Cycl_{(\xi,\eta,\zeta)}\<\xi,\d_0l^2(i^*\eta)\zeta-
    [l^2_0\eta,l^2_0\zeta]-\varpi^2_{l^2_0\eta}\zeta\>
  \end{equation}
  Using equations~\eqref{eq:0^Psi} and~\eqref{eq:dtheta0} and the fact
  that $\upsi z=z$ for all $z\in\l$, we compute the three terms in the
  right hand side of equation~\eqref{eq:CDYBEgPsi}:
  \begin{align}\label{eq:brghe1}
    \begin{split}
      \Cycl_{(\xi,\eta,\zeta)}\<\xi,\d_0l^2(i^*\eta)\zeta\>&=
      \Cycl_{(\xi,\eta,\zeta)}\<\xi,\upsi\d_0l^1(i^*\eta)\upsi^*\zeta
      +[Ai^*\eta,\upsi l^1_0\upsi^*\zeta]
      +[\upsi l^1_0\upsi^*\eta,Ai^*\zeta]\\*
      &\qquad\quad +[Ai^*\eta,Ai^*\zeta]-[A^*\eta,Ai^*\zeta]+
      \varpi^2_{Ai^*\eta}\zeta\>
    \end{split}\\*
    \label{eq:brghe2}
    \Cycl_{(\xi,\eta,\zeta)}\<\xi,[l^2_0\eta,l^2_0\zeta]\>&=
    \<\xi,[\upsi l_0^1\upsi^*\eta+Ai^*\eta-A^*\eta,
    \upsi l^1_0\upsi^*\zeta+Ai^*\zeta-A^*\zeta]\>\\*
    \label{eq:brghe3}
    \Cycl_{(\xi,\eta,\zeta)}\<\xi,\varpi^2_{l^2_0\eta}\zeta\>&=
    \<\xi,\varpi^2_{\upsi l_0\upsi^*\eta+ Ai^*\eta-\upsi
      A^*\eta}\upsi^*\zeta\>
  \end{align}

  Now, using equations~\eqref{eq:CDYBEg}, \eqref{eq:l-eqg},
  \eqref{eq:brghe1}, \eqref{eq:brghe2}, \eqref{eq:brghe3}, and the
  $\l$-equivariance of $A$, one shows that equation~\eqref{eq:CDYBEgPsi}
  reads:
  \begin{equation}
    \<\xi\otimes\eta\otimes\zeta,\varphi^2\>=
    \<\xi\otimes\eta\otimes\zeta,\upsi^{(3)}\varphi^1\>
  \end{equation}
  $\xi,\,\eta,\,\zeta\in\g_2^*$.
  Proposition~\ref{pr:gpdmorph-->qbigmorph} is thus proved.
\end{proof}

\subsection{Gauge transformations}\label{sec:gaugetransformations}
Let $\GG=U\times G\times U$ be a dynamical Poisson groupoid over
$U\subset\l^*$ with associated groupoid cocycle $P$, and let
$\Psi\colon\GG\to\GG$ be a base preserving groupoid automorphism.  Using
$\Psi$, we can transform the Poisson structure on $\GG$ into another
Poisson structure, which is dynamical if and only if $\Psi$ is
$L$-biequivariant. Such a transformation on the dynamical Poisson
structures of $\GG$ is called \emph{gauge transformation}.

Using proposition~\ref{pr:dyn Poisson groupoid morphism}, the
transformation of $P$ by $\Psi$ is given by:
\begin{equation}\label{eq:P^Psi}
  P^\Psi_{(p,x,q)}=-l_p^\Psi+\pi^\Psi_x+\Ad_xl^\Psi_q\Ad^*_x
\end{equation}
where
\begin{align}\label{eq:l^Psi}
  l_p^\Psi&=\Ad_{\sigma_p}\upsi l_p\upsi^*\Ad^*_{\sigma_p}+\theta^\Psi_p
  +\upsi\pi_{\psi^{-1}\sigma_p}\upsi^*\\*
  \pi^\Psi_x&=\upsi\pi_{\psi^{-1}x}\upsi^*
\end{align}
The map $l^\Psi$ satisfies equation~\eqref{eq:CDYBEg} for the cocycle
$\varpi^\Psi=\varpi^\upsi$ defined by equation~\eqref{eq:varpipsi}, and
for the $3$-tensor $\varphi^\Psi=\varphi^\upsi$ defined by
equation~\eqref{eq:varphipsi}. Thus,
\begin{equation}
  \Dynl(U,\G)^\Psi\subset\Dynl(U,\G^\upsi)
\end{equation}
for all $L$-biequivariant groupoid morphism $\Psi$.

\subsection{Gauge group and moduli space of dynamical $\ell$-matrices}
The group $\Map_0(U,G)^\l$ of $\l$-equivariant smooth maps $\sigma\colon
U\to G$ satisfying $\sigma_0=1$, with pointwise multiplication acts on
$\Dynl(U,\G)$ by setting
\begin{equation}\label{eq:actionsurdynl}
  l^\sigma_p=\Ad_{\sigma_p}l_p\Ad^*_{\sigma_p}+\theta^\sigma_p
  +\pi_{\sigma_p}
\end{equation}
for all $\sigma\in\Map_0(U,G)^\l$, where $\theta^\sigma$ is defined as
the right hand side of equation~\eqref{eq:theta}. It can be checked
directly that this action is a left action:
\begin{equation}
  (l^\sigma)^{\sigma'}=l^{\sigma'\sigma}
\end{equation}
This fact can also be proved at a groupoid level, by considering the
action of $\Map_0(U,G)^\l$ on the groupoid cocycles.

The subgroup $\Map^{(2)}_0(U,G)^\l$ of $\Map_0(U,G)^\l$ consisting of
smooth maps $\sigma\in\Map_0(U,G)$ such that $\T_0\sigma=0$ acts on
$\Dynl_0(U,\G)$.

There is an equivalent at the formal level: the group $\Map_0(\DD,G)^\l$
of $\l$-equivariant formal maps $\sigma\colon\DD\to G$ satisfying
$\sigma_0=1$, with (formal) pointwise multiplication acts on the space
of formal dynamical $\ell$-matrices on $\DD$ by the same
formula~\eqref{eq:actionsurdynl} (here, use Van Est's formula (\cite{G})
to give a meaning to the term $\pi_{\sigma_p}$). While the subgroup
$\Map^{(2)}_0(\DD,G)^\l$ of $\Map_0(\DD,G)^\l$ consisting of formal maps
$\sigma\in\Map_0(\DD,G)$ such that $\T_0\sigma=0$ acts on
$\Dynl_0(\DD,\G)$.

In the terminology of~\cite{EV}, the group $\Map_0(U,G)^\l$ (resp.,
$\Map_0(\DD,G)^\l$) is called the \emph{gauge group} of dynamical
$\ell$-matrices (resp., formal dynamical $\ell$-matrices).

Set $\M(\G,\l)=\Dynl(\DD,\G)/\Map_0(\DD,\G)^\l$. In the terminology
of~\cite{X}, the space $\M(\G,\l)$ is called the \emph{moduli space of
  formal dynamical $\ell$-matrices associated with the Lie
  quasi-bialgebra $\G$}.

\subsection{Reductive decomposition and Lie quasi-bialgebra compatible
  with a reductive decomposition} The statements of this section are
generalisations of those of Etingof and Schiffmann in~\cite{ES}.
Eventhough the proofs are analogous, we shall repeat the short arguments
in the body of the text for the convenience of the reader. The proof of
theorem~\ref{th:coincide=equiv} below, which is identical to that
of~\cite{ES} is reproduced in an appendix.

Recall that a \emph{reductive decomposition of the Lie algebra
  $(\g,[\,,\,])$} is a vector space decomposition $\g=\l\oplus\m$ such
that $\l$ is a subalgebra of $\g$ and $[\l,\m]\subset\m$. Let
$\G=(\g,[\,,\,],\varpi,\varphi)$ be a Lie quasi-bialgebra, $\l$ a Lie
subalgebra of $\g$ and let $U\subset\l^*$ be an $\Ad^*_L$-invariant
subset of $\l^*$ containing $0$. If $\Dynl_0(U,\G)$ is to be non empty,
then $\varpi$ must vanish on $\l$. Let $\GG=U\times G\times U$ be the
dynamical Poisson groupoid associated with some $l\in\Dynl_0(U,\G)$. It
then follows from equation~\eqref{eq:bracket_on_gstar} that
$\g^\star_0=\l\oplus\l^\perp$ is a reductive decomposition. Therefore,
if we want the Poisson groupoid dual to $\GG$ to be dynamical too, then
$\g$ must admit a reductive decomposition $\g=\l\oplus\m$. It is natural
to ask the Lie quasi-bialgebra to be compatible in some sense with this
additional structure on $\g$. This is the purpose of
definition~\ref{def:Lie q-bialg compatible} below. From now on we make
use of the following notation: we denote by $\p_\l\colon\g\to\l$ the
projection of $\g$ on $\l$ along $\m$, and by
$s=(\p_\l)^*\colon\l^*\to\g^*$ its dual. Notice that the image of $s$ is
$\m^\perp$ and that $s$ is $\l$-equivariant.

We recall the equivariant Poincar\'e lemma (for a proof, see
\eg~\cite{ES}).
\begin{lem}[equivariant Poincar\'e lemma]
  Let $\l$ be a finite dimensional Lie algebra, $V$ a finite-dimensional
  $\l$-module, and $\omega$ an $\l$-equivariant closed $k$-form on $\DD$
  with values in $V$ for some $k\geq1$. Then, there exists an
  $\l$-equivariant $(k-1)$-form $\nu$ with values in $V$ such that
  $\d^{Rh}\nu=\omega$, where $\d^{Rh}$ is the de~Rham differential
  operator of differential forms with values in $V$.
\end{lem}

\begin{prop}\label{pr:gloubiboulga}
  Let $\G=(\g,[\,,\,],\varpi,\varphi)$ be a Lie quasi-bialgebra such
  that $\varpi_\l=0$, and assume that $\g$ admits a reductive
  decomposition $\g=\l\oplus\m$. Then every formal dynamical
  $\ell$-matrix $l\in\Dynl(\DD,\G)$ is gauge-equivalent to a dynamical
  $\ell$-matrix $l'\in\Dynl(\DD,\G)$ such that $l'_0\m^\perp=0$ and
  $l'_0\l^\perp\subset\m$.
\end{prop}
\begin{proof}
  Define $\Sigma\colon\l^*\to\g$ by $\Sigma_\alpha=-\frac12\p_\l
  l_0s\alpha-\p_\m l_0s\alpha$, and set $\sigma_p=\e^{\Sigma_p}$. Since
  $l_0$ is $\l$-equivariant (because $\varpi_\l=0$), since
  $\g=\l\oplus\m$ is a reductive decomposition, and since $\Sigma_0=0$,
  $\sigma$ lies in $\Map_0(\DD,G)^\l$, and thus defines a gauge
  transformation on $\Dynl(\DD,\G)$. Now, equations~\eqref{eq:l^Psi}
  and~\eqref{eq:0^Psi} imply:
  \begin{align}
    l^\sigma_0s\alpha=l_0s\alpha-\p_\m l_0s\alpha-\p_\l l_0s\alpha=0\\*
    \p_\l l^\sigma_0\xi=\p_\l l_0\xi-\p_\l l_0\xi=0
  \end{align}
  for all $\alpha\in\l^*$ and $\xi\in\l^\perp$.
\end{proof}

\begin{thm}\label{th:coincide=equiv}
  Let $\G=(\g,[\,,\,],\varpi,\varphi)$ be a Lie quasi-bialgebra, assume
  that $\g=\l\oplus\m$ is a reductive decomposition of the Lie algebra
  $\g$, and that $\varpi_\l=0$. Let $l$ and $l'$ be two formal dynamical
  $\ell$-matrices in $\Dynl(\DD,\G)$ such that $l_0=l'_0$. Then there
  exists a gauge transformation which transforms $l$ into $l'$.
\end{thm}
\begin{proof}
  The proof is that of~\cite{ES}. See appendix~\ref{ap:proofof}.
\end{proof}

As a corollary to theorem~\ref{th:coincide=equiv}, we have:
\begin{cor}
  If the Lie algebra $\g$ possesses a reductive decomposition
  $\g=\l\oplus\m$ and if $\varpi_\l=0$ and $\varphi\equiv0\mod\l$ then
  the quotient space $\Dynl_0(U,\G)/\Map_0^{(2)}(U,G)^\l$ consists of at
  most one point.
\end{cor}

Consider the following algebraic variety
\begin{equation}
  \M_{\G,\l,\m}=\ensemble{\t\in(\wedge^2\m)^\l}
  {\varphi^\t\equiv0\mod\l}
\end{equation}
Here, $(\wedge^2\m)^\l$ denotes the set of $\l$-equivariant elements
$\t$ of $\A(\g^*,\g)$ satisfying $\t\m^\perp=0$ and
$\t\l^\perp\subset\m$. It is immediate from
equation~\eqref{eq:actionsurdynl} that if $l$ and $l'$ are gauge
equivalent, and if $l_0\m^\perp=0$, $l_0\l^\perp\subset\m$,
$l'_0\m^\perp=0$, and $l'_0\l^\perp\subset\m$, then $l_0=l'_0$.

As a corollary to proposition~\ref{pr:gloubiboulga} and
theorem~\ref{th:coincide=equiv}, we have:
\begin{cor}\label{cor:embedding}
  Let $\g=\l\oplus\m$ be a reductive decomposition of the Lie algebra,
  and let $\G=(\g,[\,,\,],\varpi,\varphi)$ be a Lie quasi-bialgebra such
  that $\varpi_\l=0$. Then, there is a well-defined embedding
  \begin{equation}
    \M(\G,\l)\longrightarrow\M_{\G,\l,\m}
  \end{equation}
  which sends a class $\mathcal{C}$ to $l_0$, where $l\in\mathcal{C}$ is
  any representative such that $l_0\m^\perp=0$ and
  $l_0\l^\perp\subset\m$.
\end{cor}

\begin{prop}
  Let $\G=(\g,[\,,\,],\varpi,\varphi)$ be a Lie quasi-bialgebra, assume
  that $\g=\l\oplus\m$ is a reductive decomposition of the Lie algebra
  $\g$, that $\<\m^\perp,\varpi_\bullet\m^\perp\>=0$, $\varpi_\l=0$, and
  that
  $\varphi\in\Alt(\l\otimes\l\otimes\l\,\oplus\,\l\otimes\m\otimes\m)$.
  Then any formal dynamical $\ell$-matrix $l\in\Dynl_0(\DD,\G)$ is gauge
  equivalent to a formal dynamical $\ell$-matrix $l'\in\Dynl_0(\DD,\G)$
  such that $l'\m^\perp\subset\l$ and $l'\l^\perp\subset\m$.
\end{prop}
\begin{proof}
  This proposition can be proved using theorems~\ref{th:coincide=equiv}
  and~\ref{th:l-matrix}, but we give here a direct proof. Let $k\geq1$
  and assume that $l\m^\perp\subset\l$ and $l\l^\perp\subset\m$ modulo
  terms of degree $\geq k$, \ie~that $\<\xi,ls\alpha\>=0$ modulo terms
  of degree $\geq k$ for all $\alpha\in\l^*$ and $\xi\in\l^\perp$.

  If $\Sigma\colon\DD\to\g$ is an $\l$-equivariant homogeneous map of
  degree $k+1$, set $\sigma=\e^\Sigma$. Then $\sigma$ is
  $\l$-equivariant too, and $l^\sigma=l+\d\Sigma i^*-(\d\Sigma)^*$
  modulo terms of degree $\geq k+1$. Let $\nu$ be the homogeneous term
  of degree $k$ of $-\p_\m ls$. Then $\nu$ is an $\l$-equivariant
  $1$-form on $\DD$ with values in $\g$, since $\varpi_\l=0$, and since
  $\g=\l\oplus\m$ is a reductive decomposition. We now show that $\nu$
  is a closed $1$-form: let $\alpha,\,\beta\in\l^*$. Then, by
  equation~\eqref{eq:CDYBEg}, and by the assumption that $\p_\m
  ls\alpha=0$ modulo terms of degree $\geq k$, we obtain:
  \begin{equation}\label{eq:nuclosed}
    \d\nu(\alpha)\beta=-[\p_\m\d l(\alpha)\beta]_{k-1}=
    \p_\m[-\d l(\beta)\alpha
    +\<s\alpha,\varpi_{l_p\bullet}s\beta\>-
    \<s\alpha\otimes s\beta\otimes1,\varphi\>]_{k-1}
  \end{equation}
  Now, by the assumptions $\<\m^\perp,\varpi_\bullet\m^\perp\>=0$ and
  $\varphi\in\Alt(\l\otimes\l\otimes\l\;\oplus\;\l\otimes\m\otimes\m)$,
  equation~\eqref{eq:nuclosed} reads
  $\d\nu(\alpha)\beta=\d\nu(\beta)\alpha$. Thus $\nu$ is an
  $\l$-equivariant closed $1$-form on $\DD$ with values in $\g$, and
  hence, by the equivariant Poincar\'e lemma, there exists a homogeneous
  $\l$-equivariant map $\Sigma\colon\DD\to\g$ of degree $k+1$ such that
  $\d\Sigma=\nu$. Now, it is easy to see that setting $\sigma=\e^\Sigma$
  yields $\<\xi,l^\sigma s\alpha\>=0$ modulo terms of degree $\geq k+1$
  for all $\alpha\in\l^*$ and $\xi\in\l^\perp$.

  The proof follows by induction as in the proof of
  theorem~\ref{th:coincide=equiv} reproduced in
  appendix~\ref{ap:proofof}.
\end{proof}

This proposition motivates the following
\begin{defn}\label{def:Lie q-bialg compatible}
  Let $\g=\l\oplus\m$ be a reductive decomposition of the Lie algebra
  $\g$, and let $\G=(\g,[\,,\,],\varpi,\varphi)$ be a Lie
  quasi-bialgebra.
  \begin{enumerate}
  \item A dynamical $\ell$-matrix $l$ satisfying $l\m^\perp\subset\l$
    and $l\l^\perp\subset\m$ is said to be \emph{compatible with the
      reductive decomposition $\g=\l\oplus\m$}.
  \item\label{itm:q-big comp} The Lie quasi-bialgebra $\G$ is said to be
    \emph{compatible with the reductive decomposition $\g=\l\oplus\m$}
    if the following three conditions hold:
    \begin{align}
      \varpi_\l&=0\\*
      \<\m^\perp,\varpi_\bullet \m^\perp\>&=0\\*
      \varphi&\in\Alt(\l\otimes\l\otimes\l\,\oplus\,
      \l\otimes\m\otimes\m\,\oplus\,\m\otimes\m\otimes\m)
    \end{align}
  \item A Lie quasi-bialgebra $\G=(\g,[\,,\,],\varpi,\varphi)$ is said
    to be \emph{canonically compatible with the reductive decomposition
      $\g=\l\oplus\m$} if $\G$ is compatible with the reductive
    decomposition $\g=\l\oplus\m$ and if $\varphi\equiv0\mod\l$.
  \item Assume that $\g$ is compatible with the reductive decomposition
    $\g=\l\oplus\m$. Let $\g'$ be a Lie algebra and $\g'=\l\oplus\m'$ a
    reductive decomposition of $\g'$, and let
    $\G'=(\g',[\,,\,]',\varpi',\varphi')$ be a Lie quasi-bialgebra
    compatible with the reductive decomposition $\g'=\l\oplus\m'$.
    \emph{A morphism of Lie quasi-bialgebra compatible with a reductive
      decomposition on $\l$ from $\G$ to $\G'$} is a Lie quasi-bialgebra
    morphism $\upsi\colon\G\to\G'$ such that $\upsi(z)=z$ for all
    $z\in\l$ and $\upsi\m\subset\m'$.
  \end{enumerate}
\end{defn}

Let $\G=(\g,[\,,\,],\varpi,\varphi)$ be a Lie quasi-bialgebra compatible
with the reductive decomposition $\g=\l\oplus\m$. First, notice that for
all $\t\in\M_{(\g,\l,\m)}$, the Lie quasi-bialgebra $\G^\t$ is
canonically compatible with the reductive decomposition $\g=\l\oplus\m$,
and that $\g^\star_0=\l\oplus\l^\perp$ is a lagrangian subalgebra of the
canonical double $\dlie^\t$ of the twisted Lie quasi-bialgebra $\G^\t$.
Also notice that the quadruple $(\l,[\,,\,]_\l,0,\p_\l^{(3)}\varphi)$ is
a cocommutative Lie quasi-bialgebra.

\section{Alekseev--Meinrenken dynamical $r$-matrix associated with a
  cocommutative Lie quasi-bialgebra}
\label{sec:cocom LQB}
In this section we study the special case $\l=\g$. We give a dynamical
$r$-matrix associated with a cocommutative Lie quasi-bialgebra in a
neighbourhood of $0\in\g^*$, which generalizes the Alekseev--Meinrenken
dynamical $r$-matrix $\rAM$ discovered in~\cite{AM00} for the compact
case, and then adapted to the quadratic case (see \eg \cite{ES,FP}). It
was first observed in~\cite{EE} that $\rAM$ only depends on the
associator. We give a simplified proof of this fact, obtained as a
corollary to the theorems of~\cite{ES,FP}.

If $E$ is a vector space, and $f\in\L(E,E)$, we denote by
$\Spectre_E(f)\subset\Cset$ the spectrum of the endomorphism $f$ of $E$.

Consider the meromorphic function on $\Cset$:
\begin{equation}\label{eq:coth...}
  F(z)=\coth(z)-\frac1z
\end{equation}
Observe that $F$ is analytic around $0$. We recall the
Alekseev--Meinrenken theorem (for a proof, see~\cite{AM00} for the
compact case, and \cite{ES,FP} for the general case):
\begin{thm}[Alekseev--Meinrenken]\label{th:AM}
  Let $(\g,[\,,\,])$ be a Lie algebra, and $\Omega\in(\S^2\g)^\g$ a
  non-degenerate invariant symmetric $2$-tensor. Then the meromorphic
  map
  \begin{equation}
    \appli{R}{\g}{\A(\g,\g)}{p}{F(\ad_p)}
  \end{equation}
  satisfies the identity:
  \begin{multline}\label{eq:CDYBE sur g}
    \d_pR(X)Y-\d_pR(Y)X-\Omega^\sharp\big(X,\d_pR(\cdot)Y\big)\\*
    -[R_p(X),R_p(Y)]+R_p([R_p(X),Y]+[X,R_p(Y)])=[X,Y]
  \end{multline}
  for all $X,\,Y\in\g$ and $p$ such that
  $\Spectre_\g(\ad_p)\cap\im\pi\Zset^*=\emptyset$.
\end{thm}

In equation~\eqref{eq:CDYBE sur g}, the term
$\Omega^\sharp\big(X,\d_pR(\cdot)Y\big)$ is the element $u$ of $\g$ such
that $(u,v)=(X,\d_pR(v)Y\big)$ for all $v\in\g$, where $(\,,\,)$ is the
non-degenerate invariant symmetric bilinear form on $\g$ associated with
$\Omega$.

As a corollary, we have:
\begin{cor}\label{cor:canonical_r-matr}
  Let $\G=(\g,[\,,\,],0,\varphi)$ be a cocommutative Lie quasi-bialgebra
  with double $\dlie$. We consider the $\Ad^*_G$-invariant open subset
  of $\g^*$ containing $0$:
  \begin{equation}
    U=\ensemble{p\in\g^*}
    {\Spectre_\dlie(\ad^\dlie_p)\cap\im\pi\Zset^*=\emptyset}
  \end{equation}
  Then, the map:
  \begin{equation}
    \appli{\rAM}{U}{\A(\g^*,\g)}{p}{F\left(\ad^\dlie_p\right)}
  \end{equation}
  is a dynamical $r$-matrix over $U$ associated with the Lie
  quasi-bialgebra $\G$, where $\ad^\dlie$ is the adjoint action of the Lie
  algebra $\dlie$.
\end{cor}

\begin{proof}
  Since $\ad^\dlie_p\g\subset\g^*$ and $\ad^\dlie_p\g^*\subset\g$ for
  all $p\in\g^*$, and since $F$ is an odd function, it is clear that
  $F(\ad^\dlie_p)\g^*\subset\g$. Also, by the invariance of the
  (canonical) bilinear form on $\dlie$, and the oddness of the function
  $F$, we have $F(\ad^\dlie_p)i_{\g^*}\in\A(\g^*,\g)$, so that $\rAM$ is
  well-defined.

  Now, if $g\in G$ and $p\in U$, using equation~\eqref{eq:AdD}, one
  finds:
  \begin{equation*}
    \Ad_g\rAM_p\Ad^*_g=\Ad^\dlie_gF\left(\ad^\dlie_p\right)\Ad^*_g=
    F\left(\ad^\dlie_{\Ad^\dlie_gp}\right)\Ad^\dlie_g\Ad^*_g
    =\rAM_{\Ad^*_{g^{-1}}p}
  \end{equation*}
  so that $\rAM$ is $\g$-equivariant.

  It only remains to show that $\rAM$ satisfies the classical dynamical
  Yang--Baxter equation. To do so, apply theorem~\ref{th:AM} to the Lie
  algebra $(\dlie,[\,,\,]_\dlie)$ with its canonical invariant
  non-degenerate symmetric bilinear form.
\end{proof}

\begin{cor}
  Let $\G=(\g,[\,,\,],\varpi,\varphi)$ be a Lie quasi-bialgebra, and
  $\DD$ the formal neighborhood of $0\in\g^*$. Then
  $\Dynl(\DD,\G)\not=\emptyset$ if and only if $\varpi$ is an exact
  $1$-cocycle. In this case, the moduli space $\M(\G,\g)$ consists of
  one point.
\end{cor}

\section{Dynamical $\ell$-matrix compatible with a reductive
  decomposition}
\label{sec:com l-matrix}
In this section, we present the most important result of this article,
namely, we give an analytic formula for an $\ell$-matrix associated with
a Lie quasi-bialgebra canonically compatible with a reductive
decomposition of a Lie algebra. As a corollary, we show that the
embedding of the moduli space of Corollary~\ref{cor:embedding} is an
isomorphism.

\begin{thm}\label{th:l-matrix}
  Let $\g=\l\oplus\m$ be a reductive decomposition of the Lie algebra
  $(\g,[\,,\,])$, and let $\G=(\g,[\,,\,],\varpi,\varphi)$ be a Lie
  quasi-bialgebra canonically compatible with the reductive
  decomposition $\g=\l\oplus\m$ with canonical double $\dlie$. Let $U$
  be the $\Ad^*_L$-invariant open subset of $\l^*$ containing $0$
  defined by:
  \begin{equation}
    U=\ensemble{p\in\l^*}
    {\Spectre_{\l\oplus\l^*}(\ad_p^{\l\oplus\l^*})
      \cap\im\pi\Zset^*=\emptyset
      \text{ and }0\not\in\Spectre_\g(\p_\g\Ad^\dlie_{\e^{-sp}}i_\g)}
  \end{equation}
  For any $p\in U$, define:
  \begin{equation}\label{eq:can l-matrix}
    \lcan_p(s\alpha+\xi)=\rAM_p\alpha
    -(\p_\g\Ad^\dlie_{\e^{-sp}}i_\g)^{-1}\p_\g\Ad^\dlie_{\e^{-sp}}\xi
  \end{equation}
  for $\alpha\in\l^*$ and $\xi\in\l^\perp\subset\g^*$, where $\rAM$ is
  the Alekseev--Meinrenken $r$-matrix associated with the Lie
  quasi-bialgebra $\left(\l,[\,,\,]_\l,0,\p_\l^{(3)}\varphi\right)$.
  Then $\lcan$ lies in $\Dynl_0(U,\G)$ and is compatible with the
  reductive decomposition $\g=\l\oplus\m$.
\end{thm}

In the sequel, we drop the suffix $\dlie$ of the adjoint action, since
no confusion is possible.

Before proving theorem~\ref{th:l-matrix}, we state three lemmas.
\begin{lem}\label{lm:ad sur ll*}
  For all $p\in U$ and for all $n\in\Nset$, we have the following
  inclusions:
  \begin{xalignat*}{2}
    \ad^{2n}_{sp}\m^\perp&\subset \m^\perp &
    \ad^{2n}_{sp}\l&\subset\l\\*
    \ad^{2n+1}_{sp}\m^\perp&\subset\l &
    \ad^{2n+1}_{sp}s\l&\subset \m^\perp\\*
    \ad^{n}_{sp}(\m\oplus\l^\perp)&\subset\m\oplus\l^\perp
  \end{xalignat*}
  In particular, $\lcan_p$, as defined by formula~\eqref{eq:can
    l-matrix}, satisfies $\lcan_p\m^\perp\subset\l$ and
  $\lcan\l^\perp\subset\m$ for all $p\in U$.
\end{lem}
\begin{proof}
  The proof is straightforward using the assumptions made on $\varpi$
  and $\varphi$.
\end{proof}

\begin{lem}\label{lm:lmatrice antisym}
  $\lcan_p$, as defined by formula~\eqref{eq:can l-matrix}, is
  skew-symmetric for all $p\in U$.
\end{lem}
\begin{proof}
  Applying lemma~\ref{lm:linalg lemma} to the vector space decomposition
  $\dlie=\g\oplus\g^*$ and to the automorphism $\Ad_{\e^{-sp}}$ of
  $\dlie$ yields:
  \begin{equation}
    \lcan_p\xi=\p_\g\Ad_{\e^{sp}}(\p_{\g^*}\Ad_{\e^{sp}}i_{\g^*})^{-1}\xi
  \end{equation}
  for all $\xi\in\l^\perp$. Since $\Ad^*_{\e^{sp}}=\Ad_{\e^{-sp}}$ and
  since $\lcan_p\xi\in\m$ by lemma~\ref{lm:ad sur ll*}, we have
  $\lcan_p\xi=-(\lcan_p)^*\xi$. Obviously, if $\alpha\in\l^*$, then
  $(\lcan_p)^*s\alpha=-\lcan_ps\alpha$. Lemma~\ref{lm:lmatrice antisym}
  is thus proved.
\end{proof}

\begin{lem}\label{lem:ad+twist}
  For all $p\in U$ and $\xi\in\l^\perp$,
  \begin{equation}\label{eq:ad+twist}
    \Ad_{\e^{-sp}}(\lcan_p\xi+\xi)\in\g^*
  \end{equation}
  More precisely, one has:
  \begin{equation}\label{eq:=ad+twist}
    \Ad_{\e^{-sp}}(\lcan_p\xi+\xi)=
    (\p_{\g^*}\Ad_{\e^{sp}}i_{\g^*})^{-1}\xi
  \end{equation}
  for all $\xi\in\l^\perp$.
\end{lem}
\begin{proof}
  Write $\Ad_{\e^{-sp}}X=\p_\g\Ad_{\e^{-sp}}X+\p_{\g^*}\Ad_{\e^{-sp}}X$
  for all $X\in\dlie$. Then,
  \begin{align*}
    \Ad_{\e^{-sp}}(\lcan_p\xi+\xi)&=-\p_\g\Ad_{\e^{-sp}}\xi
    -\p_{\g^*}\Ad_{\e^{-sp}}(\p_\g\Ad_{\e^{-sp}}i_\g)^{-1}
    \p_\g\Ad_{\e^{-sp}}\xi
    +\Ad_{\e^{-sp}}\xi\\*
    &=\p_{\g^*}\Ad_{\e^{-sp}}\xi-\p_{\g^*}
    \Ad_{\e^{-sp}}(\p_\g\Ad_{\e^{-sp}}i_\g)^{-1}
    \p_\g\Ad_{\e^{-sp}}\xi\\*
    &=\p_{\g^*}\Ad_{\e^{-sp}}(\lcan_p\xi+\xi)\\*
    \intertext{Now, by lemma~\ref{lm:lmatrice antisym}, one has:}
    \Ad_{\e^{-sp}}(\lcan_p\xi+\xi)&=
    \p_{\g^*}\Ad_{\e^{-sp}}(\p_\g\Ad_{\e^{sp}}
    (\p_{\g^*}\Ad_{\e^{sp}}i_{\g^*})^{-1}\xi+\xi)\\*
    &=\p_{\g^*}\Ad_{\e^{-sp}}\Ad_{\e^{sp}}
    (\p_{\g^*}\Ad_{\e^{sp}}i_{\g^*})^{-1}\xi\\*
    &=(\p_{\g^*}\Ad_{\e^{sp}}i_{\g^*})^{-1}\xi
  \end{align*}
  Lemma~\ref{lem:ad+twist} is thus proved.
\end{proof}

\begin{proof}[Proof of theorem~\ref{th:l-matrix}]
  We write $l$ in place of $\lcan$.

  We show that the map $l$ is $\l$-equivariant: let $h\in L$. First,
  notice that $s\Ad^*_{h^{-1}}p=\Ad_hsp$ and that
  $\Ad_h\p_\g=\p_\g\Ad_h$, since $\pi_h=0$ (see
  equation~\eqref{eq:AdD}). Clearly, if $\alpha\in\l^*$, then
  $l_{\Ad^*_{h^{-1}}p}s\alpha=\Ad_hl_p\Ad^*_hs\alpha$. If
  $\xi\in\l^\perp$, then
  \begin{align*}
    l_{\Ad^*_{h^{-1}}p}\xi&=l_{\Ad_hp}\xi\\*
    &=-(\p_\g\Ad_{h\e^{-sp}h^{-1}}i_\g)^{-1}
    \p_\g\Ad_{h\e^{-sp}h^{-1}}\xi\\*
    &=-(\p_\g\Ad_h\Ad_{\e^{-sp}}\Ad_{h^{-1}}i_\g)^{-1}
    \p_\g\Ad_h\Ad_{\e^{-sp}}\Ad_{h^{-1}}\xi\\*
    &=-\Ad_h(\p_\g\Ad_{\e^{-sp}}i_\g)^{-1}\p_\g\Ad_{\e^{-sp}}\Ad^*_h\xi\\*
    &=\Ad_hl_p\Ad^*_h\xi
  \end{align*}
  Thus, $l$ is $\l$-equivariant.

  Clearly, $l_0=0$, and $l$ is compatible with the reductive
  decomposition $\g=\l\oplus\m$. Thus, it only remains to show that $l$
  satisfies the generalization of the classical dynamical Yang--Baxter
  equation~\eqref{eq:CDYBEg alt}. Notice that we can write this equation
  as:
  \begin{multline}
    \label{eq:CDYBEg alt q-big}
    \d_pl(i^*\xi)\eta-\d_pl(i^*\eta)\xi-i\d_p\<\xi,l_\bullet\eta\>
    -[l_p\xi,l_p\eta]+l_p\p_{\g^*}[l_p\xi,\eta]
    +l_p\p_{\g^*}[\xi,l_p\eta]\\*
    -\p_\g[l_p\xi,\eta]-\p_\g[\xi,l_p\eta]
    +l_p\p_{\g^*}[\xi,\eta]=\p_\g[\xi,\eta]
  \end{multline}
  where the bracket is that of the Lie algebra $\dlie$. We now prove
  that equation~\eqref{eq:CDYBEg alt q-big} is satisfied for all
  $\xi,\,\eta\in\g^*$:
  \begin{enumerate}
  \item If $\xi,\,\eta\in \m^\perp$, equation~\eqref{eq:CDYBEg alt
      q-big} holds, by corollary~\ref{cor:canonical_r-matr} and because
    $l$ is compatible with the reductive decomposition $\g=\l\oplus\m$.
  \item\label{item:lm} If $\xi=s\alpha$ for some $\alpha\in\l^*$ and
    $\eta\in\l^\perp$, equation~\eqref{eq:CDYBEg alt q-big} becomes:
    \begin{multline}\label{eq:CDYBEg sur l m}
      \d_pl(\alpha)\eta-[l_ps\alpha,l_p\eta]+l_p[l_ps\alpha,\eta]
      +l_p\p_{\g^*}[s\alpha,l_p\eta]\\*
      -\p_\g[s\alpha,l_p\eta]
      +l_p\p_{\g^*}[s\alpha,\eta]=\p_\g[s\alpha,\eta]
    \end{multline}
    Now, using the differential of the exponential map (see
    section~\ref{sec:diff_exp} in the appendix), the fact that
    $\Ad_{\e^{-sp}}$ is an automorphism of the Lie algebra $\dlie$, and
    the property $\p_\g[\l,\g^*]=0$, we obtain:
    \begin{enumerate}[(\ref{item:lm}--a)]
    \item
      $\d_pl(\alpha)\eta=\left(\p_\g\Ad_{\e^{-sp}}i_\g\right)^{-1}\p_\g
      \left(\left[\frac{\sinh\ad_{sp}}{\ad_{sp}}s\alpha,
          \p_{\g^*}\Ad_{\e^{-sp}} (l_p\eta+\eta)\right]\right)$
      \medskip
    \item $\left[l_ps\alpha,l_p\eta\right]=-
      \left(\p_\g\Ad_{\e^{-sp}}i_\g\right)^{-1}\p_\g\left(\vphantom{\Big(}
        \left[\Ad_{\e^{-sp}}l_ps\alpha,\p_\g\Ad_{\e^{-sp}}\eta\right]+
      \right.\\*$ $\hfills\left.\left[\left(\cosh\ad_{sp}-
            \frac{\sinh\ad_{sp}}{\ad_{sp}}\right)s\alpha,
          \p_{\g^*}\Ad_{\e^{-sp}}l_p\eta\right]\right)$
      \medskip
    \item
      $l_p[l_ps\alpha,\eta]=
      \left(\p_\g\Ad_{\e^{-sp}}i_\g\right)^{-1}\p_\g\Big(-
        \left[\Ad_{\e^{-sp}}l_ps\alpha,\p_\g\Ad_{\e^{-sp}}\eta\right]+
      \\*$
      $\hfills\left.
        \left[\left(\cosh\ad_{sp}-
            \frac{\sinh\ad_{sp}}{\ad_{sp}}\right)s\alpha,
          \p_{\g^*}\Ad_{\e^{-sp}}\eta\right]\right)$
      \medskip
    \item
      $l_p\p_{\g^*}[s\alpha,l_p\eta]=
      \p_\g\left[s\alpha,l_p\eta\right]+
      \left(\p_\g\Ad_{\e^{-sp}}i_\g\right)^{-1}\p_\g\big(
      \left[\Ad_{\e^{-sp}}s\alpha,\p_\g\Ad_{\e^{-sp}}\eta\right]-\\*$
      $\hfills
      \left[\cosh\ad_{sp}s\alpha,
        \p_{\g^*}\Ad_{\e^{-sp}}l_p\eta\right]\big)$
      \medskip
    \item
      $l_p\p_{\g^*}[s\alpha,\eta]=\p_\g[s\alpha,\eta]-
      \left(\p_\g\Ad_{\e^{-sp}}i_\g\right)^{-1}\p_\g\Ad_{\e^{-sp}}
      [s\alpha,\eta]$
    \end{enumerate}
    \medskip Assembling these terms shows that equation~\eqref{eq:CDYBEg
      sur l m} is satisfied.
  \item\label{item:mm} If $\xi,\,\eta\in\l^\perp$,
    equation~\eqref{eq:CDYBEg alt q-big} becomes:
    \begin{multline}\label{eq:CDYBEg sur m m}
      -i\d_p\<\xi,l_\bullet\eta\>
      -\left[l_p\xi,l_p\eta\right]+l_p\p_{\g^*}\left[l_p\xi,\eta\right]+
      \\*
      l_p\p_{\g^*}\left[\xi,l_p\eta\right]
      -\p_\g\left[l_p\xi,\eta\right]-\p_\g\left[\xi,l_p\eta\right]
      +l_p\p_{\g^*}[\xi,\eta]=\p_\g[\xi,\eta]
    \end{multline}
    By the previous item, we know that the projection on $\l$ of this
    equation is satisfied. So it only remains to consider the projection
    on $\m$, namely:
    \begin{equation}\label{eq:CDYBEg sur m m m}
      -\left[l_p\xi,l_p\eta\right]+l_p\p_{\g^*}\left[l_p\xi,\eta\right]+
      l_p\p_{\g^*}\left[\xi,l_p\eta\right]
      -\p_\g\left[l_p\xi,\eta\right]-\p_\g\left[\xi,l_p\eta\right]
      +l_p\p_{\g^*}[\xi,\eta]\in\l
    \end{equation}
    As in the preceding item, one obtains:
    \begin{enumerate}[(\ref{item:mm}--a)]
    \item
      $[l_p\xi,l_p\eta]=\left(\p_\g\Ad_{\e^{-sp}}i_\g\right)^{-1}\p_\g
      [\p_\g\Ad_{\e^{-sp}}\xi-\p_{\g^*}\Ad_{\e^{-sp}}l_p\xi,
      \p_\g\Ad_{\e^{-sp}}\eta-\p_{\g^*}\Ad_{\e^{-sp}}l_p\eta]$
      \medskip
    \item $l_p\p_{\g^*}[l_p\xi,\eta]=
      \p_\g[l_p\xi,\eta]+\left(\p_\g\Ad_{\e^{-sp}}i_\g\right)^{-1}\p_\g
      \big(
      \left[\p_\g\Ad_{\e^{-sp}}\xi,\p_\g\Ad_{\e^{-sp}}\eta\right]+\\*$
      $\hfills
        \left[\p_\g\Ad_{\e^{-sp}}\xi,\p_{\g^*}\Ad_{\e^{-sp}}\eta\right]
        -\left[\p_{\g^*}\Ad_{\e^{-sp}}l_p\xi,
          \Ad_{\e^{-sp}}\eta\right]\big)$
      \medskip
    \item
      $l_p\p_{\g^*}[\xi,l_p\eta]=\p_\g[\xi,l_p\eta]+
      \left(\p_\g\Ad_{\e^{-sp}}i_\g\right)^{-1}\p_\g\big(
        \left[\p_\g\Ad_{\e^{-sp}}\xi,\p_\g\Ad_{\e^{-sp}}\eta\right]+\\*$
        $\hfills
        \left[\p_{\g^*}\Ad_{\e^{-sp}}l_p\xi,\p_\g\Ad_{\e^{-sp}}\eta\right]
        -\left[\Ad_{\e^{-sp}}\xi,\p_{\g^*}
          \Ad_{\e^{-sp}}l_p\eta\right]\big)$
      \medskip
    \item
      $l_p\p_{\g^*}[\xi,\eta]=\p_\g[\xi,\eta]-
      \left(\p_\g\Ad_{\e^{-sp}}i_\g\right)^{-1}\p_\g
      \left[\Ad_{\e^{-sp}}\xi,\Ad_{\e^{-sp}}\eta\right]$
    \end{enumerate}
    \medskip Using $\p_\g[\l^\perp,\l^\perp]\subset\l$, then implies
    that~\eqref{eq:CDYBEg sur m m m} holds.
  \end{enumerate}
  Thus $l$ is a dynamical $\ell$-matrix, and theorem~\ref{th:l-matrix}
  is proved.
\end{proof}

\begin{rem}
  Theorem~\ref{th:l-matrix} improves the result of~\cite{ES} as it
  provides, on the one hand (see below), explicit analytic dynamical
  $r$-matrices in each formal gauge orbit of~\cite{ES} and, on the
  other, explicit $\ell$-matrices for certain classes of Lie
  quasi-bialgebras which are not necessarily quasi-triangular nor
  cocommutative.
\end{rem}

\begin{defn}\label{df:comp l-matrix}
  Let $\G=(\g,[\,,\,],\varpi,\varphi)$ be a Lie quasi-bialgebra
  canonically compatible with a reductive decomposition $\g=\l\oplus\m$
  of the Lie algebra $\g$.
  \begin{itemize}
  \item The dynamical $\ell$-matrix provided by
    theorem~\ref{th:l-matrix} is called \emph{the canonical dynamical
      $\ell$-matrix on $\l$ associated with the Lie quasi-bialgebra $\G$
      compatible with the reductive decomposition $\g=\l\oplus\m$}, and
    will be denoted by $\lcanG{\G,\l,\m}$, or simply $\lcan$ when no
    confusion is possible.
  \item The dynamical Poisson groupoid $\GG=U\times G\times U$
    associated with the canonical $\ell$-matrix $\lcanG{\G,\l,\m}$ on
    $U$ is called the \emph{canonical dynamical Poisson groupoid
      associated with the Lie quasi-bialgebra $\G$ compatible with the
      reductive decomposition $\g=\l\oplus\m$}, and will be denoted by
    $(\GG,\G,\l,\m)^\can$.
  \end{itemize}
\end{defn}

As a corollary to theorem~\ref{th:l-matrix}, we obtain:
\begin{cor}\label{cor:isomorphism}
  Let $\g=\l\oplus\m$ be a reductive decomposition of the Lie algebra
  $\g$, and let $\G=(\g,[\,,\,],\varpi,\varphi)$ be a quasi-Lie
  bialgebra compatible with the reductive decomposition of $\g$. Then
  the embedding $\M(\G,\l)\to\M_{\G,\l,\m}$ of
  corollary~\ref{cor:embedding} is an isomorphism.
\end{cor}
\begin{proof}
  We only need to show that this map is onto. Let
  $\rho\in\M_{\G,\l,\m}$. It is easy to check that the Lie
  quasi-bialgebra $\G^\rho$ is still compatible with the reductive
  decomposition of $\g$, and since $\varphi^\rho\equiv0\mod\l$,
  theorem~\ref{th:l-matrix} provides a dynamical $\ell$-matrix
  $\lcan\in\Dynl_0(\DD,\G^\rho)$. Now, using proposition~\ref{pr:Dynl
    twist}, $l'=\lcan+\rho$ lies in $\Dynl(\DD,\G)$ and is mapped to
  $\rho$.
\end{proof}

\begin{exmp}[$\lcan$ for the cocommutative case of~\cite{ES},
  \cite{EE}, see also~\cite{AM03,FGP}]
\label{ex:cocom}
  \begin{enumerate}[(a)]
  \item\label{itm:exmpa} Let $\g=\l\oplus\m$ with $[\l,\m]\subset\m$,
    and consider the Lie quasi-bialgebra $\G=(\g,[\,,\,],0,\varphi)$
    with
    $\varphi\in\Alt(\l\otimes\l\otimes\l\,\oplus\,\l\otimes\m\otimes\m)$.
    We have
    \begin{align}
      \p_\g(\Ad_{\e^{-sp}}\xi)&=
      -\sinh\ad_{sp}\xi\in\m,\quad\forall\xi\in\l^\perp\\*
      \p_\g(\Ad_{\e^{-sp}}u)&=\cosh\ad_{sp}u\in\m,\quad\forall u\in\m
    \end{align}
    Therefore,
    \begin{equation}
      \lcan_p (s\alpha+\xi)=
      \left(\coth\ad_{sp}-\frac1{\ad_{sp}}\right)
      s\alpha+\tanh\ad_{sp}\xi
    \end{equation}
  \item Note that~(\ref{itm:exmpa}) applies to the following classical
    situation: assume that the quadratic Lie algebra $(\g,B_\g)$ admits
    a splitting
  \begin{equation}
    \g=\g_+\oplus\g_-,\quad
    [\g_\pm,\g_\pm]\subset\g_+,\quad
    [\g_+,\g_-]\subset\g_-
  \end{equation}
  such that either
  \begin{align}
    \label{eqal:(1)}
    B_\g(\g_+,\g_-)&=0
    \intertext{or}
    \label{eqal:(2)}
    \quad B_\g(\g_\pm,\g_\pm)&=0
  \end{align}
  holds. Set $\l=\g_+$, and $\m=\g_-$. Then the associator
  $\varphi=\<\Omega,\Omega\>$, where $\Omega\in(S^2\g)^\g$ is the
  invariant element associated with $B_\g$, satisfies the condition
  of~(\ref{itm:exmpa}). Note that, in either cases, the adjoint action
  of $\dlie$ appearing in the expression of $\lcan$ coincides with that
  of $\g$.

  The Cartan decomposition of a real simple Lie algebra where $\g_\pm$
  are the $\pm1$-eigenspaces of the Cartan involution (see~\cite{H})
  provides examples of type~\eqref{eqal:(1)}. On the other hand, the
  Cartan decomposition of a complex simple Lie algebra
  $\g=\ulie\oplus\im\ulie$, where $\ulie\subset\g$ is the compact real
  form of $\g$ provides an example of type~\eqref{eqal:(1)} when
  equipped with the quadratic form $\Re B_\g$, and an example of the
  lagrangian type~\eqref{eqal:(2)} when equipped with $\Im B_\g$, where
  $B_\g$ is the Killing form of $\g$. (Note that the latter case already
  appeared in the appendix of~\cite{ES}.)
  \end{enumerate}
\end{exmp}

We end this section with two propositions:
\begin{prop}\label{pr:lcan vs lcan-}
  Let $\g=\l\oplus\m$ be a reductive decomposition of the Lie algebra
  $\g$, let $\G=(\g,[\,,\,],\varpi,\varphi)$ be a Lie quasi-bialgebra
  canonically compatible with the reductive decomposition
  $\g=\l\oplus\m$. Then,
  \begin{equation}
    \lcanG{\G^-,\l,\m}_p=-\lcanG{\G,\l,\m}_{-p}
  \end{equation}
  for all $p\in U$.
\end{prop}
\begin{proof}
  Clearly, for all $\alpha\in\l^*$,
  $\lcanG{\G^-,\l,\m}_ps\alpha=\lcanG{\G,\l,\m}_{p}s\alpha=
  -\lcanG{\G,\l,\m}_{-p}s\alpha$ (since the map $F$ defined by
  equation~\eqref{eq:coth...} is skew-symmetric). Now, let
  $\xi\in\l^\perp$. Then, using the fact that the map $J$, as defined in
  section~\ref{sec:LQBOFA}, is an isomorphism, we have:
  \begin{align*}
    \lcanG{\G^-,\l,\m}_p\xi&=
    -(\p_\g\Ad^-_{\e^{-sp}}i_\g)^{-1}\p_\g\Ad^-_{\e^{-sp}}\xi\\*
    &=(\p_\g\Ad^-_{\e^{-sp}}Ji_\g)^{-1}\p_\g\Ad^-_{\e^{-sp}}J\xi\\*
    &=(\p_\g J\Ad_{\e^{-Jsp}}i_\g)^{-1}\p_\g J\Ad_{\e^{-Jsp}}\xi\\*
    &=(\p_\g\Ad_{\e^{sp}}i_\g)^{-1}\p_\g\Ad_{\e^{sp}}\xi\\*
    &=-\lcanG{\G,\l,\m}_{-p}\xi
  \end{align*}
  Proposition~\ref{pr:lcan vs lcan-} is thus proved.
\end{proof}

\begin{prop}\label{pr:ad1ad2upsietoutletsointsoin}
  Let $\g_1=\l\oplus\m_1$ and $\g_2=\l\oplus\m_2$ be reductive
  decompositions of the Lie algebras $(\g_1,[\,,\,]_1)$ and
  $(\g_2,[\,,\,]_2)$, and let $\G_1=(\g_1,[\,,\,]_1,\varpi^1,\varphi^1)$
  and $\G_2=(\g_2,[\,,\,]_2,\varpi^2,\varphi^2)$ be two Lie
  quasi-bialgebra structures on $\g_1$ and $\g_2$ canonically compatible
  with the reductive decompositions $\g_1=\l\oplus\m_1$ and
  $\g_2=\l\oplus\m_2$. Let $\upsi\colon\G_1\to\G_2$ be a morphism of Lie
  quasi-bialgebras compatible with a reductive decomposition on $\l$.
  Then we have:
  \begin{equation}
    \upsi\lcanG{\G_1,\l,\m_1}\upsi^*=\lcanG{\G_2,\l,\m_2}
  \end{equation}
\end{prop}
\begin{proof}
  For $j=1,\,2$, we denote by $s_j\colon\l^*\to\g_i^*$, the dual of the
  projection of $\g_j$ on $\l$, along $\m_j$. Since $\upsi z=z$ for all
  $z\in\l$ and since $\upsi\m_1\subset\m_2$, then $\upsi^*s_2p=s_1p$ for
  all $p\in U$. Now, using lemma~\ref{lm:ad1 ad2}, we obtain:
  $\upsi\p_{\g_1}\Ad^1_{\e^{-s_1p}}i_{\g_1}=
  \p_{\g_2}\Ad^2_{\e^{-s_2p}}\upsi i_{\g_1}$. Thus,
  \begin{equation}\label{eq:propad1ad2upsietoutletsointsoin1}
    \upsi(\p_{\g_1}\Ad^1_{\e^{-s_1p}}i_{\g_1})^{-1}=
    (\p_{\g_2}\Ad^2_{\e^{-s_2p}}i_{\g_2})^{-1}\upsi
  \end{equation}
  By using lemma~\ref{lm:ad1 ad2} again, we obtain for all
  $\xi\in\l^\perp\subset\g_2^*$:
  \begin{equation}\label{eq:propad1ad2upsietoutletsointsoin2}
    \upsi\p_{\g_1}\Ad^1_{\e^{-s_1p}}\upsi^*\xi=
    \p_{\g_2}\Ad^2_{\e^{-s_2p}}\xi
  \end{equation}
  Assembling equations~\eqref{eq:propad1ad2upsietoutletsointsoin1}
  and~\eqref{eq:propad1ad2upsietoutletsointsoin2} proves
  proposition~\ref{pr:ad1ad2upsietoutletsointsoin}.
\end{proof}

\section{Trivialization and duality}
\label{sec:duality}
\subsection{Trivial Lie algebroids}
Let $(\g,[\,,\,]_\g)$ be a Lie algebra and $M$ a manifold. Recall
(see~\cite{M}) that the trivial Lie algebroid on $M$ with vertex algebra
$\g$ is the vector bundle $A=\T M\times\g$ over $M$, where the anchor is
the projection of $A$ on $\T M$ along $\g$ and the bracket is defined as
follows: let $\sigma$ and $\sigma'$ be two sections of the vector bundle
$A$, say $\sigma=(X,x)$ and $\sigma'=(X',x')$ where $X$ and $X'$ are two
vector fields on $M$ and $x,\,x'\colon M\to\g$, and set
\begin{equation}\label{eq:trivial algebds bracket}
  [\sigma,\sigma']_A=\bigl([X,X'],X\cdot x'-X'\cdot x+[x,x']_\g\bigr)
\end{equation}
The bracket in the first component of the right hand side of
equation~\eqref{eq:trivial algebds bracket} is the bracket of vector
fields on $M$, and $X\cdot x'$ denotes the derivative of $x'$ in the
direction of $X$.

\subsection{Duality for Lie bialgebroids}
Recall that a Lie bialgebroid on a base $M$ is a pair $(A,A')$ of
algebroids $A$ and $A'$ over $M$, together with a non-degenerate pairing
between $A$ and $A'$, and a supplementary compatibility condition
(see~\cite{MX} for the explicit statement). Now let $(\GG,\{\,,\,\})$ be
a Poisson groupoid, denote by $A(\GG)$ its associated Lie algebroid, and
by $N_\GG$ the conormal bundle of the unit in $\GG$. Note that $N_\GG$
is canonically isomorphic to $A(\GG)^*$, the dual vector bundle of
$A(\GG)$. It follows from Weinstein's coisotropic calculus
(see~\cite{W}) that $N_\GG$ carries a Lie algebroid structure induced by
the Poisson bracket $\{\,,\,\}$ on $\GG$ such that the pair
$(A(\GG),N_\GG)$ is a Lie bialgebroid.

Now, let $(A,A')$ be a Lie bialgebroid over a base $M$ and denote by $a$
and by $a'$ the anchors of $A$ and $A'$ respectively. One can show
(see~\cite{MX}) the following assertions:
\begin{enumerate}
\item the pair $(A',A)$ is a Lie bialgebroid,
\item the map $a'\circ a^*$ from $\T^*M$ to $\T M$ defines a Poisson
  bivector on $M$,
\item\label{itm:changesigne} $a\circ(a')^*=-a'\circ a^*$.
\end{enumerate}
By definition (see~\cite{M2}), the dual Lie bialgebroid of the Lie
bialgebroid $(A,A')$ is the Lie bialgebroid $(A',-A)$, where $-A$ is the
Lie algebroid obtained by changing the sign of both anchor and bracket
of $A$. The sign ``$-$'' appears in the duality to keep the same induced
Poisson structure on the base $M$, as justified by
point~\eqref{itm:changesigne} above. By definition, a Poisson groupoid
dual to a Poisson groupoid $(\GG,\{\,,\,\})$ is any (connected,
source-simply-connected) Poisson groupoid $(\GG^\star,\{\,,\,\}^\star)$
such that the Lie bialgebroid $(A(\GG^\star),N_{\GG^\star})$ is the Lie
bialgebroid dual to the Lie bialgebroid $(A(\GG),N_\GG)$. The dual is
unique up to isomorphism, but may not exist (globally) in general.

\subsection{Trivialization}
Let $G$ be a connected, simply connected Lie group with Lie algebra
$\g$, $L$ a Lie subgroup of $G$ with Lie algebra $\l$ and $U$ an
$L$-invariant open subset in $\l^*$. Consider the trivial groupoid
$\GG=U\times G\times U$, with Poisson structure given by a dynamical
$\ell$-matrix $l$ associated with a Lie quasi-bialgebra
$\G=(\g,[\,,\,],\varpi,\varphi)$. Recall (\cite{LP}, see
also~\cite{BYKS} for the case $\varpi=0$) that the Lie algebroid of the
Poisson groupoid dual to $\GG$ is the vector bundle
$N(U)=U\times\l\times\g^*$ over $U$ together with the bracket on its
sections:
\begin{equation}
\begin{aligned}
  \big[(z,\xi),(z',\xi')\big]^{N(U)}_p=
  &\Big(\d_pz'\big(a^{N(U)}_p(z_p,\xi_p)\big)
    -\d_pz\big(a^{N(U)}_p(z'_p,\xi'_p)\big)-[z_p,z'_p]+
    \<\xi,\d_pl(\cdot)\xi'\>,\\*
  &\qquad
    \d_p\xi'\big(a^{N(U)}_p(z_p,\xi_p)\big)
    -\d_p\xi\big(a^{N(U)}_p(z'_p,\xi'_p)\big)
  +\ad^*_{z_p}\xi'_p-\ad^*_{z_p'}\xi_p\\*
  &\qquad+\<\xi_p,\varpi_\bullet\xi'_p\>
  +\ad^*_{l_p\xi_p}\xi'_p-\ad^*_{l_p\xi'_p}\xi_p\Big)
\end{aligned}
\end{equation}
and the anchor:
\begin{equation}
  a^{N(U)}_p(z,\xi)=i^*\xi-\ad^*_zp
\end{equation}
Since the anchor is a submersion onto $U$, a theorem of Mackenzie
(see~\cite{M}) shows that for a contractible base $U$ this algebroid is
trivializable (that is, it is isomorphic to the trivial Lie algebroid
$U\times\h^*\times\Ker a^{N(U)}_p$ for any $p\in U$).

In this section, we give an explicit trivialization for the dynamical
$\ell$-matrix of theorem~\ref{th:l-matrix} (note that, here, $U$ is not
in general contractible). First of all, we need a Lie algebra
isomorphism $\phi_p\colon\g^\star_p\to\g^\star_0$, given by the
following proposition:
\begin{prop}\label{pr:LAB iso}
  The hypotheses and notations are those of theorem~\ref{th:l-matrix}.
  For all $p\in U$, the map:
  \begin{equation}
    \appli{\phi_p}{\g^\star_p}{\g^\star_0}{X}
    {\Ad^\dlie_{\e^{-sp}}\tau_{l_p}X}
  \end{equation}
  is a Lie algebra isomorphism. Thus, the bundle map
  \begin{equation}
    \appli{\psi}{U\times\g^\star_0}
    {\Ker a^{N(U)}\subset U\times(\l\oplus\g^*)}
    {(p,X)}{\left(p,-\phi_p^{-1}X\right)}
  \end{equation}
  is a Lie algebra bundle isomorphism.
\end{prop}
As in section~\ref{sec:com l-matrix}, we drop the suffix $\dlie$ in the
adjoint actions.

\begin{proof}
  Let $p\in U$. For any $z\in\l$, since $s\ad^*_zp=\ad_{sp}z$, a
  computation yields (see lemmas~\ref{lm:ad sur ll*},
  and~\ref{lem:ad+twist}):
  \begin{equation}
    \Ad_{e^{-sp}}\tau_{l_p}(z+s\ad^*_zp)=
    \frac{\ad_{sp}}{\sinh\ad_{sp}}z\in\l
  \end{equation}
  and for any $\xi\in\l^\perp$,
  \begin{equation}
    \Ad_{e^{-sp}}\tau_{l_p}\xi=
    \left(\p_{\g^*}\Ad_{\e^{sp}}i_{\g^*}\right)^{-1}\xi\in\l^\perp
  \end{equation}
  thus $\phi_p$ is well-defined. It is clearly a Lie algebra isomorphism
  (since $\Ad_{e^{-sp}}$ and $\tau_{l_p}$ are Lie algebra isomorphisms).
\end{proof}

To complete the trivialization, we need a flat connection:
\begin{prop}\label{pr:flat connection & psi}
  Under the hypotheses and notations of theorem~\ref{th:l-matrix} and
  proposition~\ref{pr:LAB iso}, the bundle map
  \begin{equation}
    \appli{\theta}{U\times\l^*}{N(U)=U\times\l\times\g^*}
    {(p,\alpha)}
    {\left(p,\dfrac{\sinh\ad_{sp}-\ad_{sp}}{\ad_{sp}^2}s\alpha,
        \dfrac{\sinh\ad_{sp}}{\ad_{sp}}s\alpha\right)}
  \end{equation}
  is a flat connection satisfying:
  \begin{equation}\label{eq:theta compatible psi}
    \bigl[\theta\alpha,\psi X\bigr]^{N(U)}=\psi\big(\d X(\alpha)\big)
  \end{equation}
  for any smooth section $\alpha\in\Gamma(U\times\l^*)$ and
  $X\in\Gamma\left(U\times\g^\star_0\right)$.
\end{prop}
\begin{proof}
  First, notice that
  \begin{equation}\label{eq:lp sh ch}
    \frac{\sinh\ad_{sp}-\ad_{sp}}{\ad_{sp}^2}s\alpha=
    \frac{\cosh\ad_{sp}-1}{\ad_{sp}}s\alpha
    -l_p\frac{\sinh\ad_{sp}}{\ad_{sp}}s\alpha
  \end{equation}
  for all $p\in U$ and $\alpha\in\l^*$.

  For $p\in U$ and $\alpha\in\l$, an easy computation shows that
  $a^{N(U)}_p\bigl(\theta(p,\alpha)\bigr)=\alpha$. Thus, $\theta$ is a
  flat connection if and only if:
  \begin{equation}
    \bigl[\theta(\alpha),\theta(\beta)\bigr]^{N(U)}=0
  \end{equation}
  for all constant sections $\alpha,\,\beta$. For
  $\alpha,\,\beta\in\l^*$ seen as constant sections of the vector bundle
  $U\times\l^*$,
  \begin{equation}\label{eq:is theta flat?}
    \begin{aligned}
      \bigl[\theta(\alpha),\theta(\beta)\bigr]^{N(U)}_p&=
      \left(
        \d_p\frac{\sinh\ad_{s\cdot}-\ad_{s\cdot}}
        {\ad_{s\cdot}^2}(\alpha)s\beta-
        \d_p\frac{\sinh\ad_{s\cdot}-\ad_{s\cdot}}
        {\ad_{s\cdot}^2}(\beta)s\alpha\right.\\*&\qquad\left.
        -\left[
          \frac{\sinh\ad_{sp}-\ad_{sp}}{\ad_{sp}^2}s\alpha,
          \frac{\sinh\ad_{sp}-\ad_{sp}}{\ad_{sp}^2}s\beta
        \right]
        +\left<\frac{\sinh\ad_{sp}}{\ad_{sp}}s\alpha,
          \d_pl(\cdot)\frac{\sinh\ad_{sp}}{\ad_{sp}}s\beta\right>,
      \right.\\*&\qquad\left.
        \d_p\frac{\sinh\ad_{s\cdot}}{\ad_{s\cdot}}(\alpha)s\beta
        -\d_p\frac{\sinh\ad_{s\cdot}}{\ad_{s\cdot}}(\beta)s\alpha
      \right.\\*&\qquad\left.
        -\left[
          \frac{\sinh\ad_{sp}-\ad_{sp}}{\ad_{sp}^2}s\alpha,
          \frac{\sinh\ad_{sp}}{\ad_{sp}}s\beta
        \right]
        -\left[
          \frac{\sinh\ad_{sp}}{\ad_{sp}}s\alpha,
          \frac{\sinh\ad_{sp}-\ad_{sp}}{\ad_{sp}^2}s\beta
        \right]
      \right.\\*&\qquad\left.
        -\left[
          \frac{\sinh\ad_{sp}}{\ad_{sp}}s\alpha,
          \frac{\cosh\ad_{sp}}{\ad_{sp}}s\beta-
          \frac{\sinh\ad_{sp}}{\ad_{sp}^2}s\beta
        \right]
      \right.\\*&\qquad\left.
        -\left[
          \frac{\cosh\ad_{sp}}{\ad_{sp}}s\alpha-
          \frac{\sinh\ad_{sp}}{\ad_{sp}^2}s\alpha,
          \frac{\sinh\ad_{sp}}{\ad_{sp}}s\beta
        \right]
      \right)
    \end{aligned}
  \end{equation}
  The first equation of lemma~\ref{lm:diff shad/ad} shows that the
  second component of the right hand side of equation~\eqref{eq:is theta
    flat?} vanishes. To prove that the first component vanishes too, use
  equation~\eqref{eq:lp sh ch}, the generalization of the classical
  dynamical Yang--Baxter equation~\eqref{eq:CDYBEg alt},
  lemma~\ref{lm:diff shad/ad}, and the $\l$-equivariance
  equation~\ref{eq:l-eqg}. Thus $\theta$ is a flat connection.

  Now, we show that equation~\eqref{eq:theta compatible psi} is
  satisfied. As in the first part of this proof, we only need to show
  that equation~\eqref{eq:theta compatible psi} is satisfied for all
  $\alpha\in\l^*$ considered as a constant section:
  \begin{enumerate}
  \item If $X_p=z_p\in\l$ for all $p\in U$,
    \par\noindent
    $\begin{aligned}\bigl[\theta(\alpha),\psi(z)\bigr]^{N(U)}_p&=
      -\bigl[\theta(\alpha),\theta(i^*\ad_{s\cdot}z)\bigr]^{N(U)}_p
      -\bigl[\theta(\alpha),(z,0)\bigr]^{N(U)}_p\\*
      &=-\theta[\alpha,i^*\ad_{s\cdot}z]_p
      -\bigl[\theta(\alpha),(z,0)\bigr]^{N(U)}_p
      \text{, since $\theta$ is a flat connection}\\*
      &=\psi_p\bigl(\d_pz(\alpha)\bigr)\end{aligned}$
    \par\noindent
    by lemma~\ref{lm:analytic diff adsp}.
  \item A direct computation (using the fact that
    $\ad_z\p_{\g^*}=\p_{\g^*}\ad_z$ for all $z\in\l$) shows that
    equation~\eqref{eq:theta compatible psi} is satisfied for
    $X_p=\xi_p\in\l^\perp$ for all $p\in U$.
  \end{enumerate}
  Proposition~\ref{pr:flat connection & psi} is thus proved.
\end{proof}
Assembling proposition~\ref{pr:LAB iso} and proposition~\ref{pr:flat
  connection & psi}, we get:

\begin{thm}\label{th:trivialization}
  Under the hypotheses and notations of theorem~\ref{th:l-matrix}, the
  bundle map:
  \begin{equation}
    T\colon U\times\l^*\times\g^\star_0\longrightarrow
    N(U)=U\times\l\times\g^*
  \end{equation}
  given by
  \begin{equation}\label{eq:trivialisation}
    T_p(\alpha,z+\xi)=
    \left(\frac{\sinh\ad_{sp}-\ad_{sp}}{\ad_{sp}^2}s\alpha-
      \frac{\sinh\ad_{sp}}{\ad_{sp}}z,
      \frac{\sinh\ad_{sp}}{\ad_{sp}}s\alpha-\sinh\ad_{sp}z-
      \p_{\g^*}\Ad_{\e^{sp}}\xi\right)
  \end{equation}
  is a Lie algebroid isomorphism.
\end{thm}

The inverse $T^{-1}$ of $T$ is given by
\begin{equation}
  T^{-1}_p(z,s\alpha+\xi)=\left(
    \alpha-\ad^*_zp,
    \frac{\sinh\ad_{sp}-\ad_{sp}}{\ad_{sp}\sinh\ad_{sp}}s\alpha-z
    -(\p_{\g^*}\Ad_{e^{sp}}i_{\g^*})^{-1}\xi
    \right)
\end{equation}
for $z\in\l$, $\alpha\in\l^*$ and $\xi\in\l^\perp$.

\subsection{Duality}
\label{ssc:duality}
We start with a definition of a duality for Lie quasi-bialgebras:
\begin{defn}\label{df:q-big dual}
  Let $\G=(\g,[\,,\,],\varpi,\varphi)$ be a Lie quasi-bialgebra with
  canonical double $\dlie$, and assume that $\l$ is a Lie subalgebra of
  $\g$ such that $\varpi_\l=0$ and $\varphi\equiv0\mod\l$. The Lie
  quasi-bialgebra
  \begin{equation*}
    \G^\star=\left(\G_{(\dlie,\l\oplus\l^\perp,\m^\perp\oplus\m)}\right)^-
  \end{equation*}
  is called the \emph{dual over $\l$ of the Lie quasi-bialgebra $\G$}.
\end{defn}
Observe that if $\varpi_\l=0$ and $\varphi\equiv0\mod\l$ then
$\g^\star=\l\oplus\l^\perp$ is indeed a lagrangian subalgebra of
$\dlie$, so that the dual over $\l$ is well-defined. Also observe that
if a Lie quasi-bialgebra is canonically compatible with a reductive
decomposition $\g=\l\oplus\m$ (see definition~\ref{def:Lie q-bialg
  compatible}), then its dual over $\l$ is also canonically compatible
with the reductive decomposition $\g^\star=\l\oplus\l^\perp$.

Let $\op\colon\g\to\g$ be the standard involution associated to the
reductive decomposition $\g=\l\oplus\m$:
\begin{equation*}
  \op(z)=z\qquad\qquad
  \op(u)=-u
\end{equation*}
for all $z\in\l$ and $u\in\m$. Denote by $\dlie^\star$ the canonical
double of $\G^\star$. First, observe that under the canonical vector
space identification $\dlie^\star\simeq\dlie$, the Lie algebra $\g$ is
not a Lie subalgebra of $\dlie^\star$ (but the Lie algebra $\g$ is
isomorphic to the Lie subalgebra $\g^\op=\l\oplus\m$ of $\dlie^\star$).
Second, observe that under the canonical identification
$\dlie^\star\simeq\dlie$, then $(\G^\star)^\star\not=\G$, but rather
$(\G^\star)^\star=\G^\op$, which is isomorphic to $\G$.

We now turn to our main duality statement which provides the dual
Poisson groupoid of a Poisson groupoid associated with a canonical
$\ell$-matrix:
\begin{thm}\label{th:duality}
  Under the hypotheses and notations of theorem~\ref{th:l-matrix}, the
  dual Poisson groupoid of the dynamical Poisson groupoid associated
  with the canonical $\ell$-matrix $\lcan$ is (isomorphic to) the
  connected, source-simply-connected covering of the dynamical Poisson
  groupoid $U\times G^\star\times U$ with the Poisson structure
  associated with the canonical $\ell$-matrix on $U$ for the Lie
  quasi-bialgebra $\G^\star$, where $G^\star$ is the connected, simply
  connected Lie group with Lie algebra $\g^\star$.
\end{thm}
\begin{proof}
  The Poisson bracket on the dual groupoid $\GG^\star$ is uniquely
  determined (up to automorphism) by the requirement that the
  trivialization map $T$ of theorem~\ref{th:trivialization} is a Lie
  bialgebroid isomorphism, that is, by the condition that the map $-T^*$
  is a Lie algebroid isomorphism from the Lie algebroid $A(\GG)$ of
  $\GG$ to the Lie algebroid $N_{\GG^\star}(U)$ (the conormal bundle of
  the unit of the Poisson groupoid $\GG^\star$). We compute $-T^*\colon
  U\times\l^*\times\g\to U\times\l\times(\g^\star)^*\simeq
  U\times\l\times(\l^*\oplus\m)$:
  \begin{equation}\label{eq:dual trivialisation}
    -T_p^*(\alpha,z+u)=\left(
      \frac{\sinh\ad_{sp}-\ad_{sp}}{\ad_{sp}^2}s\alpha
      -\frac{\sinh\ad_{sp}}{\ad_{sp}}z,
      \frac{\sinh\ad_{sp}}{\ad_{sp}}s\alpha-\sinh\ad_{sp}z
      +\p_\g\Ad_{\e^{-sp}}u
    \right)
  \end{equation}
  where $\alpha\in\l^*$, $z\in\l$ and $u\in\m$. Here, the adjoint action
  is that of the double $\dlie$ of $\G$ on itself.

  Let $\GG^\op$ be the trivial groupoid $U\times G^\op\times U$, and
  denote by $T^\star\colon A(\GG^\op)\to U\times\l\times(\g^\star)^*$
  the trivialization associated with the data $\G^\star$ and
  $\g^\star=\l\oplus\l^\perp$ given by theorem~\ref{th:trivialization}.
  The Lie algebroid isomorphism $T^\star$ is given by:
   \begin{equation}\label{eq:trivialisation star op}
    T^\star_p(\alpha,z+u)=\left(
      \frac{\sinh\ad^\star_{sp}-\ad^\star_{sp}}{(\ad_{sp}^\star)^2}s\alpha
      -\frac{\sinh\ad^\star_{sp}}{\ad^\star_{sp}}z,
      \frac{\sinh\ad^\star_{sp}}{\ad^\star_{sp}}
      s\alpha-\sinh\ad^\star_{sp}z
      -\p_{(\g^\star)^*}\Ad^\star_{\e^{sp}}u
    \right)
  \end{equation}
  Here, $\ad^\star$ denotes the adjoint action of the double
  $\dlie^\star$ of $\G^\star$ on itself. Denote by $J^\star$ the linear
  isomorphism from $\dlie^\star$ to $\dlie$ defined by
  $J^\star(z+\xi)=z+\xi$ for all $z+\xi\in\g^\star$ and
  $J^\star(s\alpha+u)=-s\alpha-u$ for all $\alpha\in\l^*$ and $u\in\m$
  (we use the canonical vector space identification
  $\dlie^\star\simeq\dlie$). Recall that $J^\star$ is a Lie algebra
  isomorphism from $\dlie^\star$ to $\dlie$ (see
  section~\ref{sec:LQBOFA}). Using $J^\star$,
  equation~\eqref{eq:trivialisation star op} reads:
  \begin{align*}\label{eq:trivialisation star op dans dlie}
    T^\star_p(\alpha,z+u)&=\left(
      J^\star\frac{\sinh\ad^\star_{sp}-\ad^\star_{sp}}{(\ad_{sp}^\star)^2}
      s\alpha
      -J^\star\frac{\sinh\ad^\star_{sp}}{\ad^\star_{sp}}z,
      -J^\star\frac{\sinh\ad^\star_{sp}}{\ad^\star_{sp}}s\alpha
      +J^\star\sinh\ad^\star_{sp}z
      +\p_{\g}J^\star\Ad^\star_{\e^{sp}}u
    \right)\\*
    &=\left(
      \frac{\sinh\ad_{sp}-\ad_{sp}}{\ad_{sp}^2}s\alpha
      -\frac{\sinh\ad_{sp}}{\ad_{sp}}z,
      \frac{\sinh\ad_{sp}}{\ad_{sp}}s\alpha
      -\sinh\ad_{sp}z
      -\p_{\g}\Ad_{\e^{-sp}}u
    \right)
  \end{align*}
  Now, denote by $\widehat\op\colon U\times\l^*\times\g\to
  U\times\l^*\times\g^\op$ the trivial Lie algebroid isomorphism given
  by: $\widehat\op_p(\alpha,z+u)=(\alpha,z-u)$. Clearly,
  $-T^*\circ\widehat\op=T^\star$. Thus, $-T^*$ is a Lie algebroid
  isomorphism.
\end{proof}

\begin{exmp}[Dual Lie quasi-bialgebras for Etingof--Varchenko dynamical
  $r$-matrices] Let $\g$ be a complex simple Lie algebra with Killing
  form $B_\g$, $\h\subset\g$ a Cartan subalgebra, $\Delta$ (resp.,
  $\Delta^s$) the set of roots (resp., simple roots). Denote by
  $\<\Gamma\>\subset\Delta$ the root span of a fixed subset
  $\Gamma\subset\Delta^s$, and set
  $\barGamma^\pm=\Delta^\pm\setminus\<\Gamma\>^\pm$, where $\Delta^\pm$
  denotes positive and negative roots. Let
  \begin{equation}
    \g=\h\oplus\bigoplus_{\alpha\in\Delta}\g^\alpha
  \end{equation}
  be the root space decomposition of $\g$. Denote by $(x_i)_{(1\le
    i\le\rank\g)}$ an orthonormal basis of $\h$ and choose root vectors
  $(e_\alpha)_{\alpha\in\Delta}$ such that
  $B_\g(e_\alpha,e_{-\alpha})=1$.

  In~\cite{EV}, Etingof and Varchenko have shown that, up to gauging,
  analytic dynamical $r$-matrices at $0$ associated with the Lie
  quasi-bialgebra
  $\G=\left(\g,[\,,\,],0,\frac14\<\Omega,\Omega\>\right)$ are given by
  \begin{equation}
    \REV_q (\eta)=\sum_{i,j}C_{ij}(q)\<x_j,\eta\> x_i+
    \sum_{\alpha\in\Delta}\phi_\alpha(q)\<e_{-\alpha},\eta\>e_\alpha
  \end{equation}
  for all $\eta\in\g^*$, where
  \begin{align}
    \phi_\alpha(q)&=\frac12\coth\frac{(\alpha, q-\mu)}2,
    \quad\forall\alpha\in\<\Gamma\>,\\*
    \quad\phi_\alpha(q)&=\pm\frac12,
    \quad\forall\alpha\in\barGamma^\pm.
  \end{align}
  Here, $\sum_{ij}C_{ij}\d x_i\otimes\d x_j$ is an arbitrary closed
  analytic $1$-form on $\h^*$ vanishing at $0$, and $\mu\in\h^*$ lies in
  the complement of the singular hyperplanes $\<\alpha,\mu\>=0$,
  $\alpha\in\<\Gamma\>$.

  Let $\rho=\REV_0$. Note that $\REV_0$ lies in the algebraic variety
  $\M_{\G,\l,\m}$ of~\cite{ES}. By corollary~\ref{cor:isomorphism},
  $R^{EV}$ is (formally) gauge equivalent to $l^{can}+\rho$, where
  $l^{can}$ is the canonical $\ell$-matrix associated with the twisted
  Lie quasi-bialgebra
  $\mathcal{G}^\rho=(\mathfrak{g},[\,,\,],\partial\rho,\varphi^\rho)$,
  and hence, the Poisson groupoids associated with $l^{can}$ and
  $R^{EV}$ are (formally) isomorphic. In particular, the dual Poisson
  groupoid is given (up to a formal isomorphism) by
  Theorem~\ref{th:duality}.

  We now give the pair of dual Lie quasi-bialgebras
  $(\G^\rho,\G^\star)$. In the formulae below, we use the
  identifications $\g^*\simeq\g$, $\h\oplus\h^\perp\simeq\h\oplus\n$ and
  $\m^\perp\oplus\m\simeq\h\oplus\n$ induced by the Killing form.

  For $\G^\rho$, we have, for all $x,\,y\in\g$,
  \begin{align*}
    \dell\rho_xy&=[x,\rho(y)]_\g-\rho ([x,y]_\g)\\*
    (x\otimes y\otimes1,\varphi^\rho)&=
    \frac14[x,y]_\g+[\rho(x),\rho(y)]_\g
    -\rho\bigl([\rho(x),y]_\g+[x,\rho(y)]_\g\bigr)
  \end{align*}
  The structural data for $\G^\star=
  (\g^\star=\h\oplus\h^\perp,[\,,\,]^\star,\varpi^\star,\varphi^\star)$
  are as follows: the Lie bracket $[\,,\,]^\star$ of
  $\h\oplus\h^\perp\simeq\h\oplus\n$ is that of the double $\dlie$ for
  $\G^\rho$:
  \begin{align*}
    [z,z']^\star&=0\\*
    [z,u]^\star&=[z,u]_\g\\*
    [u,u']^\star&=(u\otimes u'\otimes1,\varphi^\rho)
    +(u',\dell\rho_\cdot u)
  \end{align*}
  for $z,\,z'\in\h$ and $u,\,u'\in\n$.

  Put $\phi_\alpha=\phi_\alpha(0)$. The bracket
  $[\,,\,]^\star_{\mid_{\n\times\n}}$ then reads as
  \begin{align*}
    [e_\alpha,e_\beta]^\star&=(\phi_\alpha+\phi_{\beta})\,
    [e_\alpha,e_\beta]_\g\\*
    [e_\alpha,e_{-\alpha}]^\star&=(\frac14-\phi_\alpha^2)\,
    [e_\alpha,e_{-\alpha}]_\g
  \end{align*}
  for all $\alpha, \beta\in \Delta$ such that $\alpha+\beta\not=0$,
  which implies that $\g^\star$ is isomorphic to the semi-direct product
  $\l_\Gamma\ltimes(\n_\Gamma^+\ominus\n_\Gamma^-)$, where
  $$\l_\Gamma=\h\oplus\bigoplus\limits_{\alpha\in\<\Gamma\>}\g^\alpha$$
  is the Levi factor and $\n_\Gamma^\pm$ are the corresponding nilpotent
  radicals (a fact already observed in~\cite{LP}).

  The cocycle $\varpi^*$ is given by
  $-\p_{\h\oplus\h^\perp}[z+u,z'+u']_\dlie$, for
  $z+u\in\h\oplus\h^\perp$, and $z'+u'\in\m^\perp\oplus\m$:
  \begin{equation}
    \varpi^*_{z+u}z'+u'=\p_\h\dell\rho_{u'}u
    -\p_\h(u\otimes z'\otimes 1,\varphi^\rho)
    -(z',\dell\rho_\cdot u)+\p_{\n}[u',u]_\g
  \end{equation}
  which explicitely reads as
  \begin{align*}
    \varpi^\star_z&=0\\*
    \varpi^\star_{e_\alpha}z&=\phi_\alpha\,[z,e_\alpha]_\g\\*
    \varpi^\star_{e_\alpha}e_\beta&=-[e_\alpha,e_\beta]_\g\\*
    \varpi^\star_{e_\alpha}e_{-\alpha}&=
    -\phi_\alpha\,[e_\alpha,e_{-\alpha}]_\g
  \end{align*}
  for $z\in\h$, and for all $\alpha,\,\beta\in\Delta$ such that
  $\alpha+\beta\not=0$.

  The associator $\varphi^\star$ is given by
  $\p_{\h\oplus\h^\perp}[z+u,z'+u']_\dlie$, for
  $z+u,\,z'+u'\in\m^\perp\oplus\m$:
  \begin{align*}
    (u\otimes u'\otimes 1,\varphi^\star)&=\p_\h[u,u']_\g\\*
    (u\otimes z'\otimes1, \varphi^\star)&=
    \p_\h\dell\rho_u z'+[u,z']_\g=[u,z']_\g\\*
    (z\otimes z'\otimes1,\varphi^\star)&=
    \p_\h(z\otimes z'\otimes1,\varphi^\rho)
    +(z',\dell\rho_\cdot z)=0
  \end{align*}

  One checks that, when $\Gamma\not=\Delta$, the cocycle
  $\varpi^\star\colon\g^*\to\L(\m^\perp\oplus\m,\g^*)$ is \emph{not}
  exact, thus providing a genuine example of non-exact Lie
  quasi-bialgebra which is compatible with a reductive decomposition.
\end{exmp}

\subsection{Link with the duality of symmetric space}
In this section, we show that the duality of symmetric spaces
(see~\cite{H}), which relies on duality for orthogonal symmetric Lie
algebras, is related to the duality of quasi-bialgebra introduced in
definition~\ref{df:q-big dual}. We first recall the definition of
orthogonal symmetric Lie algebras (see~\cite{H}).
\begin{defn}
  An \emph{orthogonal symmetric Lie algebra} is a pair $(\g,\sigma)$
  where $\g$ is a Lie algebra over $\Rset$, $\sigma$ is an involutive
  automorphism of $\g$, and $\l$, the set of fixed points of $\sigma$ is
  a compactly imbedded subalgebra of $\g$.
\end{defn}

When $(\g,\sigma)$ is an orthogonal symmetric Lie algebra, we consider
the reductive splitting $\g=\l\oplus\m$ of $\g$ into the eigenspaces of
$\sigma$ for the eigenvalue $+1$ and $-1$ respectively as in
example~\ref{ex:cocom}.

There is a notion of duality for orthogonal symmetric Lie algebras which
goes as follows: let $(\g,\sigma)$ be an orthogonal symmetric Lie
algebra. Let $\g^\star$ denote the subset $\l+\im\m$ of the
complexification $\g^\Cset$ of $\g$. Notice that $\g^\star$ is a (real)
Lie subalgebra of $\g^\Cset$, since $[\m,\m]\subset\l$. Now, the mapping
$\sigma^\star\colon z+\im u\mapsto z-\im u$ is an involutive
automorphism of $\g^\star$, and the pair $(\g^\star,\sigma^\star)$ is
again an orthogonal symmetric Lie algebra, called the \emph{dual of the
  orthogonal symmetric Lie algebra $(\g,\sigma)$}.

We can translate this duality when $\g$ is semi-simple: let
$(\g,\sigma)$ be an orthogonal symmetric Lie algebra, and assume that
$\g$ is a semi-simple Lie algebra. We consider the complexified Lie
algebra $\dlie=\g^\Cset$, which we view as a real Lie algebra. We
consider the bilinear, symmetric, invariant, and non-degenerate bilinear
form on $\dlie$: $(\,,\,)_\dlie=\Im B$, where $B$ is the Killing form on
$\g^\Cset$. Then $\g$ is a lagrangian Lie subalgebra of $\dlie$, and
$\im\g$ is an isotropic complement. We set $\G=\G_{(\dlie,\g,\im\g)}$.
Clearly, the Lie quasi-bialgebra $\G$ is cocommutative, and its
associator $\varphi$ is given by $\varphi=-\<\Omega,\Omega\>$, where
$\Omega\in(S^2\g)^\g$ is the Casimir element of $\g$, and $\<\,,\,\>$ is
Drinfel\cprime d's bracket. Now, since $\Omega$ lies in
$\l\otimes\l\oplus\m\otimes\m$, and since $[\m,\m]\subset\l$, the Lie
quasi-bialgebra $\G$ is canonically compatible with the reductive
decomposition $\g=\l\oplus\m$. The dual of the orthogonal symmetric Lie
algebra $(\g,\sigma)$ is $(\g_0^\star,\sigma^\star)$, where $\g_0^\star$
is the underlying Lie algebra of $\G^\star$, the dual of the
quasi-bialgebra $\G$, and where $\sigma^\star$ is the standard
involution associated to the reductive decomposition on $\g^\star_0$.
We can note that $\g^\star_0$ is still semi-simple (use Cartan's
criterion), but that $\g$ and $\g^\star_0$ are not isomorphic in
general: this can be seen by comparing the signature of the Killing form
on $\g$ and on $\g^\star_0$. Thus, the Lie quasi-bialgebra $\G$ is not
self-dual, nor is the dynamical Poisson groupoid $\GG$ associated to
$\G$.

\begin{appendix}
\section{Complements}
\subsection{A linear algebra result}
To prove that the map defined by equation~\eqref{eq:can l-matrix} is
skew-symmetric, we use the following lemma:
\begin{lem}\label{lm:linalg lemma}
  Let $E$ be a vector space, $F$ and $F'$ two subspaces such that
  $E=F\oplus F'$, and denote by $\p$ and $\p'$ the projections on $F$
  and $F'$ along $F'$ and $F$ and by $i$ and $i'$ the inclusions of $F$
  and $F'$ into $E$. Let $f$ be an automorphism of $E$, and assume that
  $\p fi$ and $\p'f^{-1}i'$ are automorphisms of $F$ and $F'$. Then:
  \begin{equation}\label{eq:easylemma}
    (\p fi)^{-1}\p fi'=-\p f^{-1}i'(\p'f^{-1}i')^{-1}
  \end{equation}
\end{lem}
\begin{proof}
  We start with the relation:
  \begin{equation}\label{eq:easyrelation}
    \L(F',F)\ni0=\p ff^{-1}i'=(\p fi)\p f^{-1}i'+\p fi'(\p'f^{-1}i')
  \end{equation}
  Since $\p fi$ and $\p'f^{-1}i'$ are automorphisms of $F$ and $F'$, one
  obtains equation~\eqref{eq:easylemma} by multiplying $(\p fi)^{-1}$
  and $(\p'f^{-1}i')^{-1}$ on the left and on the right of
  equation~\eqref{eq:easyrelation}.
\end{proof}

\subsection{Differential of the exponential map}\label{sec:diff_exp}
For the convenience of the reader we recall the expression of the
differential of the exponential map (for a proof see \eg\cite{H}):

Let $G$ be a Lie group with Lie algebra $\g$, and denote by
$\exp\colon\g\to G$ the exponential map. For any $x\in\g$, the
differential of $\exp$ at $x$ is given by
\begin{equation}
  \T_x\exp(u)=\ell_{\exp x}\frac{1-\e^{-\ad_x}}{\ad_x}u
\end{equation}
for all $u\in\g$. In particular, one has:
\begin{equation}
  \d_x\Ad_{\e^\bullet}(u)=\Ad_{\e^x}\ad_{\frac{1-\Ad_{\e^{-x}}}{\ad_x}u}
\end{equation}
for all $u\in\g$.

\subsection{More differential identities}
We prove some differential identities which are used to prove
proposition~\ref{pr:flat connection & psi}:
\begin{lem}\label{lm:diff shad/ad}
  For all $p\in U$ and $\alpha,\,\beta\in\l^*$, the following three
  equations hold:
  \begin{equation}
    \label{eq:diff ad^n}
    \d_x\ad^n(u)v=\sum_{i=0}^{n-1}\binom n{i+1}[\ad^i_xu,\ad^{n-i-1}_xv]
  \end{equation}
  \begin{multline}\label{eq:dshsh}
    \d_p\frac{\sinh\ad_{s\cdot}}{\ad_{s\cdot}}(\alpha)s\beta
    -\d_p\frac{\sinh\ad_{s\cdot}}{\ad_{s\cdot}}(\beta)s\alpha=\\*
    \left[\frac{\cosh\ad_{sp}-1}{\ad_{sp}}s\alpha,
      \frac{\sinh\ad_{sp}}{\ad_{sp}}s\beta\right]+ \left[
      \frac{\sinh\ad_{sp}}{\ad_{sp}}s\alpha,
      \frac{\cosh\ad_{sp}-1}{\ad_{sp}}s\beta\right]
  \end{multline}
  \begin{multline}\label{eq:dchch}
    \d_p\frac{\cosh\ad_{s\cdot}-1}{\ad_{s\cdot}}(\alpha)s\beta
    -\d_p\frac{\cosh\ad_{s\cdot}-1}{\ad_{s\cdot}}(\beta)s\alpha=\\*
    \left[\frac{\sinh\ad_{sp}}{\ad_{sp}}s\alpha,
      \frac{\sinh\ad_{sp}}{\ad_{sp}}s\beta\right]+\left[
      \frac{\cosh\ad_{sp}-1}{\ad_{sp}}s\alpha,
      \frac{\cosh\ad_{sp}-1}{\ad_{sp}}s\beta\right]
  \end{multline}
\end{lem}
\begin{proof}
  To prove equation~\eqref{eq:diff ad^n}, write:
  \begin{equation}
    \d_x\ad^n(u)v=\sum_{i=0}^{n-1}\ad_x^i[u,\ad_x^{n-i-1}v]
  \end{equation}
  and use Leibniz' relation:
  \begin{equation}
    \ad^n_x[u,v]=\sum_{i=0}^n\binom n{i}[\ad^i_xu,\ad^{n-i}_xv]
  \end{equation}
  and the identity
  \begin{equation}
    \sum_{i=0}^{n-k}\binom{i+k}k=\binom{n+1}{k+1}
  \end{equation}
  Recall that the expressions $\frac{\sinh\ad_{sp}}{\ad_{sp}}$ and
  $\frac{\cosh\ad_{sp}-1}{\ad_{sp}}$ are respectively the even and odd
  parts of $\frac{\Ad_{\e^{sp}}-1}{\ad_{sp}}$. Using lemma~\ref{lm:diff
    shad/ad}, we compute:
  \par\medskip\noindent
  $\begin{aligned}
    \d_p\frac{\Ad_{\e^{s\cdot}}-1}{\ad_{s\cdot}}(\alpha)s\beta-
    \d_p\frac{\Ad_{\e^{s\cdot}}-1}{\ad_{s\cdot}}(\beta)s\alpha
    &=\sum_{n\geq0}\frac1{(n+2)!}
    \left(\d_p\ad_{s\cdot}^{n+1}(\alpha)s\beta
      -\d_p\ad_{s\cdot}^{n+1}(\beta)s\alpha\right)\\*
    &=\sum_{n\geq0}\sum_{i=0}^n\frac1{(n+2)!}
    \binom{n+2}{i+1}[\ad_{sp}^is\alpha,\ad_{sp}^{n-i}s\beta]\\*
    &=\left[\frac{\Ad_{\e^{sp}}-1}{\ad_{sp}}s\alpha,
      \frac{\Ad_{\e^{sp}}-1}{\ad_{sp}}s\beta\right]
  \end{aligned}$

  Now, selecting respectively odd and even parts of this expression
  yields relations~\eqref{eq:dshsh} and~\eqref{eq:dchch}.
\end{proof}

The following lemma is used to prove relation~\eqref{eq:theta compatible
  psi} of proposition~\ref{pr:flat connection & psi}.
\begin{lem}\label{lm:analytic diff adsp}
  For all entire functions $f$, the identity
  \begin{equation}
    \left.\frac{\d}{\d t}\right\vert_{t=0} f(\ad_{sp+t\ad_{sp}z})
    =f(\ad_{sp})\ad_z-\ad_zf(\ad_{sp})
  \end{equation}
  holds for all $z\in\l$ and for all $p\in U$.
\end{lem}
\begin{proof}
  Let $f$ be an entire function, say $f(x)=\sum_{n\geq0}f_nx^n$. Then,
  for $v\in\dlie$ and $z\in\l$, using lemma~\ref{lm:diff shad/ad} one
  computes:
  \begin{equation*}
    \begin{aligned}
      \d_{\ad_{sp}}f(\ad_{sp}z)v&=
      \sum_{n\geq1}f_n\d_{sp}\ad^n(\ad_{sp}z)v\\*
      &=\sum_{n\geq1}f_n\sum_{i=0}^{n-1}
      \binom ni[\ad^i_{sp}z,\ad_{sp}^{n-i}v]\\*
      &=\sum_{n\geq0}f_n\left(\ad_{sp}^n[z,v]-[z,\ad_{sp}^nv]\right)
    \end{aligned}
  \end{equation*}
  Lemma~\ref{lm:analytic diff adsp} is thus proved.
\end{proof}

\section{Proof of lemma~\ref{lm:ad1 ad2}}
\label{ap:pflem}
Observe that equation~\eqref{eq:lm ad1 ad2 -- 1d} is the dual of
equation~\eqref{eq:lm ad1 ad2 -- 1}, so the two are equivalent. We prove
lemma~\ref{lm:ad1 ad2} by induction on $n$: for $n=1$, the
relations~\eqref{eq:lm ad1 ad2 -- 1}, \eqref{eq:lm ad1 ad2 -- 2}
and~\eqref{eq:lm ad1 ad2 -- 3} hold, by definition of a quasi-bialgebra
morphism. Assume that these relations hold for some $n\in\Nset$. Then,
\begin{align*}
  \begin{split}
    \upsi\p_{\g_1}\left(\ad^1_{\upsi^*\xi}\right)^{n+1}u
    &=\upsi\p_{\g_1}\ad^1_{\upsi^*\xi}\p_{\g_1}
    \left(\ad^1_{\upsi^*\xi}\right)^nu
    +\upsi\p_{\g_1}\ad^1_{\upsi^*\xi}\p_{\g_1^*}
    \left(\ad^1_{\upsi^*\xi}\right)^nu\\*
    &=\p_{\g_2}\ad_\xi^2\p_{\g_2}\Bigl(\ad_\xi^2\Bigr)^n\upsi u
    +\upsi\p_{\g_1}\ad^1_{\upsi^*\xi}\upsi^*\p_{\g_2^*}
    \Bigl(\ad_\xi^2\Bigr)^n\upsi u\\*
    &=\p_{\g_2}\ad_\xi^2\p_{\g_2}\Bigl(\ad_\xi^2\Bigr)^n\upsi u
    +\p_{\g_2}\ad_\xi^2\p_{\g_2^*}\Bigl(\ad_\xi^2\Bigr)^n\upsi u
  \end{split}
  \intertext{which proves~\eqref{eq:lm ad1 ad2 -- 1} at rank $n+1$,}
  \begin{split}
    \p_{\g_1^*}\left(\ad^1_{\upsi^*\xi}\right)^{n+1}u
    &=\p_{\g_1^*}\ad^1_{\upsi^*\xi}\p_{\g_1}
    \left(\ad^1_{\upsi^*\xi}\right)^nu
    +\p_{\g_1^*}\ad^1_{\upsi^*\xi}\p_{\g_1^*}
    \left(\ad^1_{\upsi^*\xi}\right)^nu\\*
    &=\upsi^*\p_{\g_2^*}\ad^2_\xi\upsi\p_{\g_1}
    \left(\ad^1_{\upsi^*\xi}\right)^nu
    +\p_{\g_1^*}\ad^1_{\upsi^*\xi}\upsi^*\p_{\g_2^*}
    \Bigl(\ad_\xi^2\Bigr)^n\upsi u\\*
    &=\upsi^*\p_{\g_2^*}\ad^2_\xi\p_{\g_2} \Bigl(\ad_\xi^2\Bigr)^n\upsi
    u +\upsi^*\p_{\g_2^*}\ad^2_\xi\p_{\g_2^*}
    \Bigl(\ad_\xi^2\Bigr)^n\upsi u
  \end{split}
  \intertext{which proves~\eqref{eq:lm ad1 ad2 -- 2} at rank $n+1$, and}
  \begin{split}
    \upsi\p_{\g_1}\left(\ad^1_{\upsi^*\xi}\right)^{n+1}\upsi^*\eta
    &=\upsi\p_{\g_1}\ad^1_{\upsi^*\xi}\p_{\g_1}
    \left(\ad^1_{\upsi^*\xi}\right)^n\upsi^*\eta
    +\upsi\p_{\g_1}\ad^1_{\upsi^*\xi}\p_{\g_1^*}
    \left(\ad^1_{\upsi^*\xi}\right)^n\upsi^*\eta\\*
    &=\p_{\g_2}\ad^2_\xi\upsi\p_{\g_1}
    \left(\ad^1_{\upsi^*\xi}\right)^n\upsi^*\eta
    +\upsi\p_{\g_1}\ad^1_{\upsi^*\xi}\upsi^*\p_{\g_2^*}
    \Bigl(\ad_\xi^2\Bigr)^n\eta\\*
    &=\p_{\g_2}\ad^2_\xi\p_{\g_2}\Bigl(\ad_\xi^2\Bigr)^n\eta
    +\p_{\g_2}\ad^2_\xi\p_{\g_2^*}\Bigl(\ad_\xi^2\Bigr)^n\eta
  \end{split}
\end{align*}
which proves~\eqref{eq:lm ad1 ad2 -- 3} at rank $n+1$.
Lemma~\ref{lm:ad1 ad2} is thus proved.

\section{Proof of theorem~\ref{th:coincide=equiv}}
\label{ap:proofof}
Let $k\geq1$, and assume that $l=l'$ modulo terms of degree $\geq k$.
We show that there exists $\sigma\in\Map_0(\DD,G)^\l$ such that
$l^\sigma=l'$ modulo terms of degree $\geq k+1$: for all
$\xi\in\l^\perp$, one has $\p_\m l\xi=\p_\m l'\xi$ modulo terms of
degree $\geq k+1$. Indeed, equation~\eqref{eq:CDYBEg} shows that the
term of degree $k-1$ of $\<\xi,\d l(\alpha) \eta\>$ only depends on the
terms of degree $\leq k-1$ of $l$. Thus, for all $\alpha\in\l^*$ and
$\xi,\,\eta\in\l^\perp$, one has $\<\xi,\d l(\alpha) \eta\>=\<\xi,\d
l'(\alpha) \eta\>$ modulo terms of degree $\geq k$, and since
$l_0=l'_0$, the equality $\p_\m l\xi=\p_\m l'\xi$ holds modulo terms of
degree $\geq k+1$. Now, if $\Sigma\colon\DD\to\g$ is an $\l$-equivariant
homogeneous map of degree $k+1$, set $\sigma=\e^\Sigma$. Then $\sigma$
is $\l$-equivariant, and one checks that $l^\sigma=l+\d\Sigma
i^*-(\d\Sigma)^*$ modulo terms of degree $\geq k+1$. We show that there
exists such a map $\Sigma$ such that $l^\sigma=l'$ modulo terms of
degree $\geq k+1$: we define a $2$-form $\mu$ on $\DD$ (with values in
the ground field) by setting
$\<\mu,\alpha\wedge\beta\>=[\<s\alpha,(l'-l)s\beta\>]_k$, for
$\alpha,\,\beta\in\l^*$, and we define a $1$-form $\nu$ on $\DD$ with
values in $\g$ by setting $\<\nu,\alpha\>=\p_\m[(l'-l)]_ks\alpha$, for
$\alpha\in\l^*$. Here, the bracket $[\cdot]_k$ means to select the
homogeneous term of degree $k$. Since $\varpi_\l=0$ and $\g=\l\oplus\m$
is a reductive decomposition, both $\mu$ and $\nu$ are $\l$-equivariant
forms (the scalar field is seen as a trivial $\l$-module). Using
equation~\eqref{eq:CDYBEg} and the assumption that $l=l'$ modulo terms
of degree $\geq k$, we check that they are closed forms:
\begin{align*}
  \begin{split}
    \d\<\mu,\alpha\wedge\beta\>(\gamma)&=
    [\<s\alpha,\d(l'-l)(\gamma)s\beta\>]_{k-1}\\*
    &=[\<s\alpha,\d(l'-l)(\beta)s\gamma\>]_{k-1}+
    [\<s\gamma,\d(l'-l)(\alpha)s\beta\>]_{k-1}\\*
    &=\d\<\mu,\alpha\wedge\gamma\>(\beta)
    +\d\<\mu,\gamma\wedge\beta\>(\alpha)
  \end{split}\\*
  \d\<\nu,\alpha\>(\beta)&=\p_\m[\d(l'-l)(\beta)s\alpha]_{k-1}=
  \p_\m[\d(l'-l)(\alpha)s\beta]_{k-1}\\*
  &=\d\<\nu,\beta\>(\alpha)
\end{align*}
Thus, $\mu$ is an $\l$-equivariant closed $2$-form on $\DD$ and $\nu$ is
an $\l$-equivariant closed $1$-form on $\DD$ with values in $\g$. Thus,
by the equivariant Poincar\'e lemma, there exists a homogeneous
$\l$-equivariant $1$-form $\chi$ on $\DD$ of degree $k+1$ such that
$\d^{Rh}\chi=\mu$, and a homogeneous $\l$-equivariant map
$\lambda\colon\DD\to\g$ of degree $k+1$ such that $\d\lambda=\nu$. The
$1$-form $\chi$ may be seen as an $\l$-equivariant map from $\DD$ to
$\l$, which will be denoted by $\overline\chi$. Now, set
$\Sigma_p=\overline\chi_p+\lambda_p$. We check that setting
$\sigma=\e^\Sigma$ yields $l^\sigma=l'$ modulo terms of degree $\geq
k+1$: let $\alpha,\,\beta\in\l^*$. Clearly,
$\<s\beta,\d\Sigma(\alpha)\>-\<s\alpha,\d\Sigma(\beta)\>=
\<\d^{Rh}\chi,\alpha\wedge\beta\>=[\<s\alpha,(l'-l)s\beta\>]_k$, so that
$\<s\beta,(l^\sigma-l')\alpha\>=0$ modulo terms of degree $\geq k+1$.
Let $\alpha\in\l^*$ and $\xi\in\l^\perp$. Then,
$\<\xi,\d\Sigma(\alpha)\>=\<\xi,\d\lambda(\alpha)\>=
\<\xi,[(l'-l)]_ks\alpha\>$, so that $\<\xi,(l^\sigma-l')s\alpha\>=0$
modulo terms of degree $\geq k+1$. Therefore, $l^\sigma=l'$ modulo terms
of degree $\geq k+1$.

It is clear that $l=l'$ modulo terms of degree $\geq1$. Thus, by
induction we construct a sequence $\sigma^{(k)}\in\Map_0(\DD,G)^\l$ of
homogeneous map of degree $k$ such that
$l^{\sigma^{(k)}\cdots\sigma^{(2)}}=l'$ modulo terms of degree $\geq k$.
Clearly, the sequence $\sigma^{(k)}\cdots\sigma^{(2)}$ converges in
$\Map_0(\DD,G)^\l$ to a map $\sigma\in\Map_0(\DD,G)^\l$ such that
$l^\sigma=l'$. Theorem~\ref{th:coincide=equiv} is thus proved.
\end{appendix}


\end{document}